\documentclass[12pt]{article}
\usepackage{amsfonts}
\usepackage{amsmath}
\usepackage{verbatim}
\usepackage{amssymb}
\usepackage[latin1]{inputenc}
\usepackage{a4wide}
\usepackage{amsthm}
\usepackage{graphicx}

\DeclareMathAlphabet{\mathpzc}{OT1}{pzc}{m}{it}

\newcommand{\N}{\ensuremath{\mathbb N}}

\newcommand{\R}{\ensuremath{\mathbb R}}

\newcommand{\T}{\ensuremath{\mathbb T}}
\newcommand{\Y}{\ensuremath{\mathbb Y}}
\newcommand{\Pp}{\ensuremath{\mathbb P}}

\newcounter{comptage}[part]
\newtheorem{lem}[comptage]{Lemma}
\newtheorem{theo}[comptage]{Theorem}
\newtheorem{defin}[comptage]{Definition}
\newtheorem{cor}[comptage]{Corollary}
\newtheorem{prop}[comptage]{Proposition}
\newtheorem{conj}[comptage]{Conjecture}

\newtheoremstyle{remarque}
  {3pt}
  {3pt}
  {}
  {}
  {\bf}
  {.}
  {.5em}
  {}
\theoremstyle{remarque}
\newtheorem{rem}[comptage]{Remark}

\author{\vspace{1cm} Antoine Lemenant \\
Université Paris XI \\
antoine.lemenant@math.u-psud.fr}

\title{Regularity of the singular set for Mumford-Shah minimizers in $\R^3$ near a minimal cone.}

\begin{document}
\maketitle

{\bf Abstract.} We show that if $(u,K)$ is a minimizer of the
Mumford-Shah functional in an open set $\Omega$ of $\R^3$, and if
$x\in K$ and $r>0$ are such that $K$ is close enough to a minimal
cone of type $\Pp$ (a plane), $\Y$ (three half planes meeting with
120° angles) or $\T$ (cone over a regular tetrahedron centered at
the origin) in terms of Hausdorff distance in $B(x,r)$, then $K$
is $C^{1,\alpha}$ equivalent to the minimal cone in $B(x,cr)$
where $c<1$ is an universal constant.




\section*{Introduction}

The Mumford-Shah functional comes from an image segmentation
problem. If $\Omega$ is an open subset of $\R^2$, for example a
rectangle, and $g\in L^\infty(\Omega)$ is an image, D.Mumford and
J.Shah \cite{ms} proposed to define
$$J(K,u):= \int_{\Omega \backslash K}|\nabla u|^2dx+\int_{\Omega\backslash K}(u-g)^2dx+H^{1}(K)$$
and, to get a good segmentation of the image $g$, minimize the
functional $J$ over all the admissible pairs $(u,K)\in
\mathcal{A}$ (see definitions after). Any solution $(u,K)$ that
minimizes $J$ represents a ``smoother'' version of the image and
the set $K$ represents the edges of the image.

Existence of minimizers is  a well known result (see for instance
\cite{dcl}) using $SBV$ theory.

The question of regularity for the singular set $K$ of a minimizer
is more difficult. The following conjecture from D. Mumford and J.
Shah is currently still open.

\begin{conj}[Mumford-Shah]{\rm \cite{ms} } Let $(u,K)$  be a reduced minimizer for the functional $J$. Then
$K$ is the finite union of  $C^1$ arcs.
\end{conj}

Some partial results are true for this conjecture. For instance it
is known that  $K$ is $C^1$ almost everywhere (see \cite{d1},
\cite{b} and \cite{afp1}).

 Furthermore it is known that if  $B$ is a ball such that
  $K\cap B$ is a $C^{1,\alpha}$ graph, and if in addition $g$ is of regularity $C^k$,
  then $K\cap B$ is $C^k$ (cf Theorem 7.42 in \cite{afp}) and even
  that
  if  $g$ is an analytic function, then  $K$ is also analytic (see \cite{klm}).

Many results about the Mumford-Shah functional are about $\R^2$.
In dimension $3$, lots of proprieties are still unknown. The
theorem of  L. Ambrosio, N. Fusco and D. Pallara \cite{afp1} about
regularity of minimizers is one of the rare result valid in any
dimension.  It says in particular that if $K$ is flat enough in a
ball $B$, and if the energy there is not too big, then $K$ is a
$C^1$ hypersurface in a slightly smaller ball. The proof of this
result is based on a ``tilt-estimate" and does not seem to
generalize to other geometric situations different than a
hyperplane.

It is natural to think about situation in  dimension 3. Some works
on minimal surfaces and soap bubbles in dimension 3 tell us what
can be the singularities of a Mumford-Shah minimizer, at least
when the energy is small. In particular in Jean Taylor's work
\cite{ta} we can find the description of the three minimal cones
in $\R^3$. Jean Taylor also proves that any minimal surface is
locally $C^1$ equivalent to one of those cones. So we can think
that for Mumford-Shah minimizers a similar descritption is true.

What we prove here is that if in a ball, the singular set of a
Mumford-Shah minimizer is close enough to a minimal cone, then it
is $C^{1,\alpha}$ equivalent to this cone. It is a generalization
to cones $\Y$ and $\T$ of what L. Ambrosio, N. Fusco et D. Pallara
have done with hyperplanes in \cite{afp1}. It is also a
generalization in higher dimension of what G. David \cite{d1} did
in $\R^2$ about the regularity near lines and propellers.

We start with a few definitions. Let $\Omega$ be an open set of
$\R^N$. We consider the set of admissible pairs
$$\mathcal{A}:=\{(u,K); \; K \text { closed } , \; u \in W^{1,2}_{loc}(\Omega \backslash K) \}.$$

\begin{defin} Let $(u,K) \in  \mathcal{A}$ and $B$ a ball such that $\bar B\subset \Omega$.
A competitor  for the pair  $(u,K)$ in the ball $B$ is a pair
$(v,L) \in \mathcal{A}$ such that
$$
\left.
\begin{array}{c}
u=v \\
K=L
\end{array}
\right\} \text{ in } \Omega \backslash  B
$$
and in addition such  that if $x$ and $y$ are two points in
$\Omega \backslash ( B \cup K) $ that are separated by $K$ then
they are also separated by $L$.
\end{defin}

The expression ``be separated by $K$'' means that $x$ and $y$ lie
in different connected components  of $\Omega \backslash K$.

\begin{defin}
A gauge function $h$ is a non negative and non decreasing function
 on $\R^+$ such that $\lim_{t \to 0}h(t)=0$.
\end{defin}

\begin{defin}\label{defms} Let $\Omega$ be an open set of $\R^N$.  A Mumford-Shah
minimizer with gauge function $h$ is a pair $(u,K)\in \mathcal{A}$
such that for every ball $\bar B \subset \Omega$ and every
competitor $(v,L)$ in $B$ we have
$$\int_{B \backslash K}|\nabla u|^2dx +H^{N-1}(K\cap B)\leq \int_{B\backslash L}|\nabla v|^2dx+H^{N-1}(L\cap B)+r^{N-1}h(r)$$
with $r$ the radius of the ball $B$ and where $H^{N-1}$ denotes
the Hausdorff measure of dimension $N-1$.
\end{defin}

It is not difficult to prove that a minimizer for the functional
$J$ of the beginning of the introduction is a minimizer in the
sense of Definition \ref{defms} with $h(r)=C_N\|g\|_{\infty}^2r$
as gauge function where $C_N$ is a dimensional constant (see
proposition 7.8 p. 46 of \cite{d}).

\begin{defin}\label{defmsg} A global minimizer in $\R^N$ is a
Mumford-Shah minimizer in the sense of Definition \ref{defms} with
$\Omega=\R^N$ and $h=0$.
\end{defin}

We will not work on global minimizers in this paper but they take
an important place in the study of the Mumford-Shah functional and
that is why we introduced the definition. In dimension 2, only
three types of connected sets can give a global minimizer : $K$ is
a line and $u$ is locally constant, $K$ is a propeller (a union of
three half-lines meeting with 120 degree angles) and $u$ is
locally constant as well, and finally when $K$ is a half line and
$u$ is a $cracktip$. Knowing whether there is another global
minimizer would give a positive answer to the Mumford-Shah
conjecture. The main fact is that every blow up limit of $(u,K)$
is a global minimizer. In \cite{l1}, one can find some
informations about global minimizers in $\R^3$.

If $(u,K)$ is a Mumford-Shah minimizer and if we add to $K$ a
small closed set of Hausdorff measure zero, then  this new set is
also a Mumford-Shah minimizer. That is why in the following we
will always suppose that the minimizer is ``reduced". This means
that a pair  $(\tilde u,\tilde K) \in \mathcal{A}$ such that
$\tilde K \varsubsetneq K$ and $\tilde u$ is an extension of $u$
in $W^{1,2}_{loc}(\Omega \backslash \tilde K)$ doesn't exist.
Given a pair $(u,K) \in \mathcal{A}$, one can always find a
reduced pair $(\tilde u , \tilde K) \in \mathcal{A}$ such that
$\tilde K \subset K
  $ and $\tilde u$ is an extension of $u$ (see Proposition 8.2 of \cite{d}).

Let us now define the minimal cones that will be used in the next
sections. We define three types of cones. Cones of type $1$ are
planes in $\R^3$, also called $\Pp$. Cones of types $2$ and $3$
and their spines are defined as in  \cite{ddpt} by the following
way.

\begin{defin}\label{prop} Define $Prop\subset \R^2$
by
$$Prop=\{(x_1,x_2);x_1 \geq 0, x_2=0\} $$
$$\hspace{4cm} \cup\{(x_1,x_2);x_1 \leq 0, x_2=-\sqrt{3}x_1\}$$
$$\hspace{6.5cm}\cup\{(x_1,x_2);x_1 \leq 0, x_2=\sqrt{3}x_1\}.$$
Then set $Y_0=Prop\times \R \subset \R^3.$ The spine of $Y_0$ is
the line $L_0=\{x_1=x_2=0\}$. A cone of type $2$ (or of type $\Y$)
is a set $Y=R(Y_0)$ where  $R$  is the composition of a
 translation and a rotation. The spine of $Y$ is then the line $R(L_0)$.
 We denote by $\Y$ the set of all the cones of type 2. Sometimes we also may use
 the expression ``of type $\Y$".
\end{defin}

\begin{defin}\label{defT} Let $A_1=(1,0,0)$, $A_2=(-\frac{1}{3},\frac{2\sqrt{2}}{3},0)$,
$A_3=(-\frac{1}{3},-\frac{\sqrt{2}}{3},\frac{\sqrt{6}}{3})$, and
$A_4=(-\frac{1}{3},-\frac{\sqrt{2}}{3}, -\frac{\sqrt{6}}{3})$ the
four vertices of a regular tetrahedron centered at $0$. Let $T_0$
be the cone over the union of the $6$ edges $[A_i,A_j]$ $i\not
=j$. The spine of $T_0$  is the union of the four half lines
$[0,A_j[$. A cone of type $3$ (or of type $\T$) is a set
$T=R(T_0)$ where $R$ is the composition of a translation and a
rotation. The spine of $T$ is the image by $R$  of the spine of
$T_0$. We denote by $\T$ the set of all the cones of type 3.
\end{defin}

\begin{center}
\includegraphics[width=6cm]{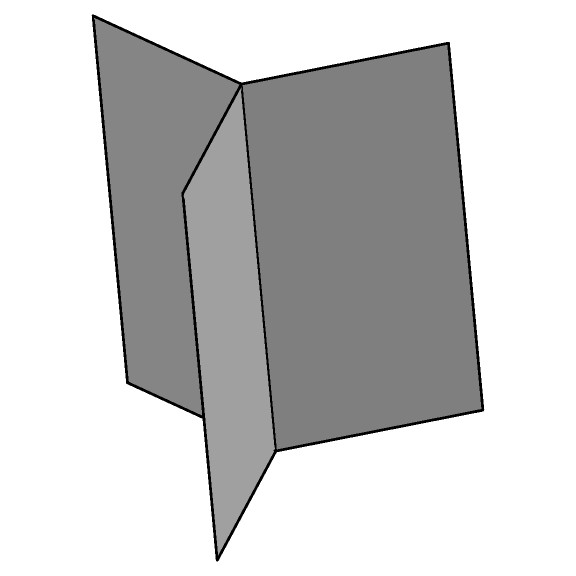} \hspace{2cm}
\includegraphics[width=6cm]{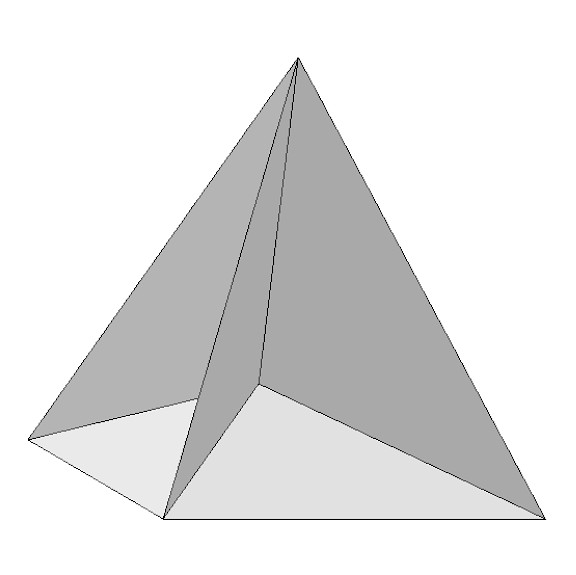}\\
Cones\footnote{Thanks to Ken Brakke for those pictures.} of type
$\Y$ and $\T$.
\end{center}

Cones of type $\Pp$, $\Y$ and $\T$ are the only sets (except the
empty set) in $\R^3$ that locally minimizes the Hausdorff measure
of dimension 2 under topological conditions (i.e. every competitor
keep the same connected components outside the competitor ball).
This fact is proved in \cite{d5}. That is why in the following we
will say ``minimal cones'' to design cones of type $\Pp$, $\Y$ and
$\T$.

We denote by $D_{x,r}$ the normalized Hausdorff distance between
two closed sets $E$ and $F$ in $B(x,r)$ defined by
$$D_{x,r}(E,F):=\frac{1}{r}\Big\{\max\{\sup_{y\in E\cap B(x,r)}d(y,F),\sup_{y\in F\cap B(x,r)}d(y,E)\}\Big\}.$$

We now come to the main result of the paper.

\begin{theo}\label{fin} We can find some absolute positive constants $\varepsilon>0$ and $c<1$ such that all
the following is true. Let $(u,K)$ be a reduced Mumford-Shah
minimizer in $\Omega \subset \R^3$, with
 gauge function $h$. Let $x \in K$ and $r>0$ be such that
$B(x,r)\subset{\Omega}$. Assume in addition that there is a
minimal cone $Z$ of type $\Pp$, $\Y$ or $\T$ centered at $x$ such
that
$$D_{x,r}(K,Z)+h(r)\leq \varepsilon.$$
 Then there is a diffeomorphism $\phi$ of class $C^{1,\alpha}$ from $B(x, cr)$ to its image such that $K \cap B(x, cr)=\phi(Z)\cap
B(x, cr)$.
\end{theo}

When $(u,K)$ is a Mumford-Shah minimizer in $\Omega\subset \R^N$,
and if $B(x,r)$ is a ball such that $\bar B(x,r)\subset \Omega$,
we denote by $\omega_2(x,r)$ the normalized energy of $u$ in
$B(x,r)$ defined by
$$\omega_{2}(x,r):=\frac{1}{r^{N-1}}\int_{B(x,r)\backslash K}|\nabla u|^2dx.$$

We also have a version of Theorem \ref{fin} with only a condition
on the normalized energy instead of the geometric condition.

\begin{theo}\label{finte} We can find some absolute positive constants $\varepsilon>0$ and $c<1$ such that
the following is true. Let $(u,K)$ be a reduced Mumford-Shah
minimizer in $\Omega \subset \R^3$, with
 gauge function $h$. Let $x \in K$ and $r>0$ be such that
$B(x,r)\subset{\Omega}$ and
$$\omega_2(x,r)+h(r)\leq \varepsilon.$$
 Then there is a diffeomorphism $\phi$ of class $C^{1,\alpha}$ from $B(x, cr)$
 to its image, and there is a minimal cone $Z$ such that $K \cap
B(x, cr)=\phi(Z)\cap B(x, cr)$.
\end{theo}

In all the following we will work in $\R^3$. However, the proof of
Theorem \ref{fin} still works in higher dimension for the case of
hyperplanes so that we could have a new proof of  L. Ambrosio, N.
Fusco, D. Pallara's entire result. With the same proof we could
also imagine to have other results in $\R^N$, but the analogue of
Jean Taylor's Theorem in higher dimension is missing at the time
when this paper is written.

Indeed, one of the ingredients of the proof of Theorem \ref{fin}
is to apply some results about minimal sets. In particular we will
use the paper of G. David \cite{d5} following J. Taylor \cite{ta},
that is the analogue of Theorem \ref{fin} but for almost minimal
sets. Let $E$ be a closed set in $\R^N$.

\begin{defin}\label{m1} A MS-competitor for the closed set $E$ in $\Omega\subset \R^N$ is a closed set $F$ such that
there is a ball $B\subset \Omega$ of radius $r$ with
$$F\backslash B=E\backslash B$$
and if $x,y \in \Omega\backslash (B\cup E)$ are separated by $E$
then they are also separated by $F$.
\end{defin}

\begin{defin}\label{minim}\label{m2} A set $E\subset \Omega$ is MS-almost minimal with gauge function $h$
if
$$H^{N-1}(E\cap B)\leq H^{N-1}(F\cap B)+r^{N-1}h(r)$$
for all MS-competitor $F$ for $E$ in the ball $B$ of radius $r$.
\end{defin}

If $E$ is a MS-almost minimal set, we denote
$$\theta(x,r)=r^{-2}H^2(E\cap B(x,r)).$$
 The limit at $0$ of $\theta$ exists because
$E$ is almost minimal so one can prove that $\theta$ is almost non
decreasing (see 2.3 of \cite{d5}). The limit is called ``density"
of $E$ at point $x$ and will be denoted by $d(x)$. Then we
introduce the excess of density defined by
$$f(x,r)=\theta(x,r)-\lim_{t\to 0}\theta(x,t)=\theta(x,r)-d(x).$$
Now Proposition 12.28 of \cite{d5} says the following.

\begin{theo} \label{gd3} For each choice of $b \in (0,1]$, $\bar c>0$
and  $C_0>0$ we can find $\varepsilon_1>0$ and $C\geq 0$ such that
the following holds. Let $E$ be a reduced MS-almost minimal set in
$\Omega \subset \R^3$ with gauge function $h$. Suppose that $0 \in
E$, $r_0>0$ be such that $B(0,110r_0)\subset \Omega$ and $h$ is
satisfying
$$h(r)\leq C_0r^b\quad \text{ for } 0<r<220r_0.$$
Assume in addition that
$$f(0,110r_0)+C_0r_0^b\leq \varepsilon_1$$
and
$$D_{0,100r_0}(E,Z)\leq \varepsilon_1$$
where $Z$ is a minimal cone centered at the origin such that
$$H^2(Z\cap B(0,1))\leq d(0).$$
Then for all $x\in E$ and $r>0$ such that $x \in E\cap B(0,10r_0)$
and $0<r<10r_0$, we can find a minimal cone $Z(x,r)$, not
necessarily centered at $x$ or at the origin, such that
$$D_{x,r}(E, Z(x,r))\leq \bar c \left(\frac{r}{r_0}\right)^{\alpha}$$
\end{theo}

The constant $\alpha$ is a universal constant depending on
dimension and other geometric facts.

We also will need this result (Corollary 12.25 of  \cite{d5}).

\begin{cor}{\rm \cite{d5}} \label{gd} For each choice of $b \in (0,1]$,
and  $C_0>0$ we can find $\alpha>0$ and $\varepsilon_1>0$ such
that the following holds. Let $E$ be a reduced MS-almost minimal
set in $\Omega \subset \R^3$ with gauge function $h$. Suppose that
$0 \in E$, $r_0>0$ is such that $B(0,110r_0)\subset \Omega$ and
$h$ is satisfying
$$h(r)\leq C_0r^b\quad \text{ for } 0<r<220r_0.$$
Assume in addition that
\begin{eqnarray}
f(0,110r_0)+C_0r_0^b\leq \varepsilon_1 \label{defautdens}
\end{eqnarray}
 and
$$D_{0,100r_0}(E,Z)\leq \varepsilon_1$$
where $Z$ is a minimal cone centered at the origin such that
$$H^2(Z\cap B(0,1))\leq d(0).$$
Then for $x \in E \cap B(0,r_0)$ and $0<r\leq r_0$ there is a
$C^{1,\alpha}$ diffeomorphism $\Phi : B(0,2r) \to \Phi(B(x,2r))$,
such that $\Phi(0)=x$, $|\Phi(y)-y-x|\leq 10^{-2}r$ for $y \in
B(0,2r)$ and $E\cap B(x,r)=\Phi(Z)\cap B(x,r)$.
\end{cor}

The strategy to obtain our main result is to control the
normalized energy of $u$ (that is the quantity $\omega_2$).  While
the energy is under control with a decay as a power of radius $r$,
we can say that our singular set is a MS-almost minimal with a
gauge function that depends on the decay of $\omega_2$, thus we
can apply Corollary \ref{gd}.

We claim that if we had some similar statements as Theorem
\ref{gd3} and Corollary \ref{gd} in higher dimension, then the
work in this paper should give a analogous result for the singular
set of a minimizer for the Mumford-Shah functional in dimension
$N>3$. Unfortunately, if Guy David is quite able to give similar
results for sets of dimension 2 in $\R^N$, the technics used to
prove Theorem \ref{gd3} and Corollary \ref{gd} seem not to work
for lower co-dimensions.

The paper is organized as follow. In a first part we explain a
method to construct a good competitor using a stoping time
argument. This construction will use some preliminary work like
the Whitney extension and geometric lemmas that are also used in
\cite{l3} and which statements are recalled here. We begin by a
good control of the normalized Jump in order to avoid some
topological and geometric problems.

In the second Section we will use the competitor described in
Section 1 in order to get some estimates about  the two main
quantities that will appear: normalized energy and  bad mass. We
also prove that the minimality defect depends on those quantities.

Finally in last section we prove the decay estimate that leads to
regularity. At the end we state a few different versions of the
main theorem.

The author wishes to thank Guy David for interesting discussions
about the subject of this work and to have given some useful
remarks and suggestions about the redaction of this paper.



\section{Construction of a competitor}

It will be convenient to work with a set that is separating. That
is why in a first part we have to control the jump of function
$u$, that will be useful to estimate the size of holes in $K$.
Before that, let us recall some definitions and geometric results
from \cite{l3}.

\begin{defin}[Almost Centered] Let $Z$ be a minimal cone and  $B$ a ball that meets  $Z$. We say that $Z$ is
almost centered with constant $V$ if the center of $Z$ lies in
$\frac{1}{V}B$. If $V=2$ we just say that $Z$ is almost centered
in $B$.
\end{defin}

This lemma will be useful to deal with almost centered cones.

\begin{lem}{\rm \cite{l3}}\label{recentrage} Let $Z$
be a minimal cone in $\R^3$ that contains $0$ (but is not
necessarily centered at $0$). Then for all $r_0>0$ and for all
constant $V\geq 1$ there is a $r_1$ such that
$$r_1\in \{r_0,Vr_0, V^2 r_0\}$$
and such that we can find a cone $Z'$, containing $0$ and centered
in $B(0,\frac{1}{V}r_1)$ with $Z\cap B(0,r_1)=Z'\cap B(0,r_1)$.
\end{lem}

\begin{defin}[Separating]\label{separex} Let $Z$ be a minimal cone in
$\R^3$ and  $B$ a ball of radius $r$ such that  $B\cap Z\not =
\emptyset$. For all $a>0$ we define $Z_{a}$ by
$$Z_{a}:= \{ y \in B ; d(y,Z)\leq a\}.$$
 Let $E$ be a closed set in $B$ such that $E$ is contained in
$Z_{r\varepsilon_0}$ for some $\varepsilon_0<10^{-5}$. We say that
``$E$ is separating in $B$'' if the connected components  of
$B\backslash Z_{r\varepsilon_0}$  are contained in different
connected components of $B\backslash E$. We denote by
$\mathpzc{k}^B$ the number of connected component of $B\backslash
Z_{2\varepsilon_0}$ (thus $\mathpzc{k}^B$ is equal to $type(Z)+1$
if $Z$ is not centered too close to $\partial B$).
\end{defin}

\subsection{Separation and control of the Jump}

So let  $(u,K)$ be a Mumford-Shah minimizer in $\Omega
\subset\R^3$ (see Definition \ref{defms}) with gauge function $h$
and let  $\varepsilon$ be fixed. Suppose that there is a ball
$B(x_0,r_0)$ such that in this ball, $K$ is $\varepsilon$-close to
a minimal cone $Z^0$ of type $\Pp$,$\Y$, or $\T$ (see Definition
\ref{prop} and \ref{defT}), in other words there is a minimal cone
$Z^0$ such that
$$K\cap B(x_0,r_0) \subset Z^{0}_{\varepsilon r_0} := \{y; d(y,Z^0)\leq \varepsilon r_0\}.$$

Equivalently we have
$$\beta(x_0,r_0)\leq \varepsilon$$
where $\beta$ is the ``generalized Peter Jones unilateral
quantity" defined by
$$\beta(x,r):= \frac{1}{r}\inf_{Z}\big\{\sup\{d(y,Z); y \in K \cap B(x,r)\}\big\}.$$
The infimum is taken over all the cones of type $\Pp$, $\Y$, or
$\T$ that contain $x_0$ (but are not necessarily centered at
$x_0$).  Sometimes we will use the notation $\beta_K(x,r)$ to
precise that the quantity is associated to the set $K$.

Moreover we suppose that $Z^0$ is centered at $x_0$. Throughout
all this part of the paper, we will always work under these above
hypothesis. We introduce now some additional notations. We denote
by $\mathpzc{k}^0$ the number of connected components of
$B(x_0,r_0)\backslash Z^0$ and for all $k \in \N \cap
[1,\mathpzc{k}^0]$ we consider a ball $D_k$ of radius
$\frac{1}{10}r_0$ such that each $D_k$ are situated in one of the
connected components of $B(x_0,r_0)\backslash Z^0$, the farthest
as possible from $Z^0$. We also denote by $m_k$ the mean value of
$u$ on $D_k$. Then we introduce
$$\delta_{k,l}(x_0,r_0)=|m_k-m_l|$$
and finally, the normalized jump is defined by
$$J(x_0,r_0):=r_0^{-\frac{1}{2}}\min\{ \delta_{k,l}; 1<k,l<\mathpzc{k}^0 \text{ and } k\not = l \}.$$

In general, for all $x\in K$ and  $r>0$ such that $B(x,r)$ is
included in $\Omega$ and such that there is a cone  $Z$ almost
centered in $B(x,t)$ and $10^{-5}$ close to $K$ in $B(x,r)$, we
can define the normalized jump by the same way
$$J(x,r):=r^{-\frac{1}{2}}\min\{ \delta_{k,l}; 1<k,l<\mathpzc{k}^{B(x,r)} \text{ and } k\not = l \}.$$

Here the  $\delta_{k,l}$ are again defined as differences between
mean values of $u$ on balls of radius equivalent to $r$ in each
connected components of $B(x,r)\backslash Z$ far from $Z$.

If a ball $B(x,r)$ is such that $\beta(x,r)\leq 10^{-5}$ but with
minimal cone that realize the infimum not almost centered, we can
also define the normalized jump. Indeed, we know by the
recentering Lemma \ref{recentrage} that $B(x,2r)$ or $B(x,4r)$ is
associated to an almost centered cone. Then we define the
normalized jump
 $J(x,r)$ as being equal to the jump of the first ball between $B(x,2r)$ or $B(x,4r)$ for which the
 cone  is
 almost centered.

All the parameters that define the jump (choice of cone $Z^0$,
constant $4$ to have the almost centering property, diameter and
position of the $D_k$) are not so important since the difference
is just multiplying the jump by a constant.

First of all, we want to work with a new set $F$ that contains $K$
and such that $F$ is separating in $B(x_0,r_0)$ (see Definition
\ref{separex}). The result is the same as Proposition 1 p. 303 of
\cite{d} but generalized to the case of $\Y$ and $\T$. We also use
the opportunity here to prove an additional fact about the set $F$
(called Property $\star$) that will be used later. Recall that the
normalized energy in the ball $B$ is denoted by
$$\omega_2(x,r):=\frac{1}{r^2}\int_{B(x,r)\backslash K}|\nabla u|^2dx.$$

\begin{prop} \label{constrF}Let $(u,K)$ be a Mumford-Shah minimizer in $\Omega \subset \R^3$.
 Suppose that there is an $x\in \Omega$, a
$r>0$ and a positive constant $\varepsilon<10^{-10}$  such that
$B(x,r)\subset \Omega$ and suppose in addition that there is a
minimal cone $Z$ almost centered in $B(x,r)$ such that
$$\sup_{y \in K\cap B(x,r)}\frac{1}{r}d(y,Z)\leq \varepsilon.$$
Moreover, assume that $J(x,r)\not = 0$,
\begin{eqnarray}
\omega_2^{\frac{1}{2}}(x,r)J^{-1}(x,r)\leq \varepsilon \label{pr4}
\end{eqnarray}
 and that
\begin{eqnarray}
\omega_2(x,r)^{\frac{1}{8}}\leq CJ(x,r) \label{pr5}
\end{eqnarray}
 with  $C$ a positive universal constant given by the demonstration.
 We call $D_k$ for $k\in \N\cap
[1,\mathpzc{k}^{B(x,r)}]$ the domains in the definition of
$J(x,r)$. Then there is a compact set $F(x,r)\subset B(x,r)$ such
that
\begin{eqnarray}
K\cap B(x,r) \subset F(x,r) \subset \{x \in B ; d(x, Z) \leq Cr \sqrt{\varepsilon} \} \label{pr1} \\
F \text{ is separating each } D_k \text{ from }  D_l \text{ for  } k\not = l \text{ in } B(x,r) \label{pr2} \\
H^{2}(F(x,r)\cap B(x,r) \backslash K) \leq C
r^2\omega_2(x,r)^{\frac{1}{2}}J(x,r)^{-1} \label{pr3}\notag
\end{eqnarray}
Moreover $F$ is satisfying Property  $\star$ (defined just after).
\end{prop}

Property $\star$ shows that we control the geometry of $F$ at
small scales when the geometry of $K$ is controlled. This is the
definition.

\begin{defin}[Property $\star$] $F$ satisfy Property $\star$ if, for all $\varepsilon_0<10^{-5}$, $y \in K\cap
B(x,r)$ and $s>0$ such that
$$ \inf \{t; \forall t'\geq t,\beta_K(y,t')\leq \varepsilon_0\}\leq s \leq d(y, \partial B(x,r))$$
we have
$$\beta_F(y,s)\leq \varepsilon_0.$$
\end{defin}

\begin{rem}\label{rembete} Condition \eqref{pr5} allows us to have
Property $\star$ and Condition \eqref{pr4} is here to prove the
last inclusion of \eqref{pr1}. Proposition \ref{constrF} is still
true without Property $\star$ and without Conditions \eqref{pr4}
and \eqref{pr5}. In this case, \eqref{pr1} is proved by the use of
a retraction as in 44.1 of \cite{d}.
\end{rem}

{\bf Proof :} The first step is the same as Proposition 1 p. 303
of \cite{d} but applied to $\Y$ and $\T$ as well. However we will
write the entire proof here because it will be easier next to show
Property $\star$.

 For all $\lambda$ we call
$$S(\lambda):=\{y \in B(x,r); d(y,Z)\leq \lambda r\}$$
and denote by $A_k(\lambda)$ for  $k \in \N\cap [1,type(Z)+1]$ the
connected component of $B(x,r)\backslash S(\lambda)$ which meets
$D_k$. We set $V=B(x,r)\backslash K$. Let us find a function $v$
such that
\begin{eqnarray}
v(y)=m_k \quad \quad \text{ for } y \in A_k\left(1/10\right)
\label{condition1}
\end{eqnarray}
and \begin{eqnarray} \int_{V}|\nabla v| \leq C \int_{V}|\nabla u|.
\label{condition2}
\end{eqnarray}
To do this we consider for all $k$  a function $\varphi_k$ such
that $0\leq \varphi_k\leq 1$ and $\varphi_k=1$ on $A_k(1/10)$,
$\phi_k=0$ on $V \backslash A_k(1/100)$ and $|\nabla \phi_k|\leq C
r^{-1}$. Then we set
$$\varphi = 1-\sum_k \varphi_k$$
and
$$v=\varphi u + \sum_k \varphi_k m_k$$
while $m_k$ is the average of $u$ on $D_k$. We have
\eqref{condition1} trivially. Concerning  \eqref{condition2} we
have
$$\nabla v(y)=\varphi(y)\nabla u(y)-\sum_k {\bf \rm 1}_{A_k(1/100)}(y)\nabla \varphi_k(y)[u(y)-m_k]$$
and since $\varepsilon<10^{-5}$, the $A_k(1/100)$ do not meet $K$
and then applying Poincaré inequality in $A_k(1/100)$ gives
\begin{eqnarray}
\int_{A_k(1/100)}|\nabla \varphi_k(y)||u(y)-m_k|dy&\leq&
Cr^{-1}\int_{A_k(1/100)}|u(y)-m_k|dy \notag \\
&\leq&C \int_{A_k(1/100)}|\nabla u(y)|dy \notag
\end{eqnarray}
and \eqref{condition2} is verified.\\
Now we want to replace  $v$ with a smooth function $w$ in $V$ such
that
\begin{eqnarray}
w(y)=m_k \quad \quad \text{ for } y \in A_k(1/10)
\label{condition1p}
\end{eqnarray}
and \begin{eqnarray} \int_{V}|\nabla w| \leq C \int_{V}|\nabla u|.
\label{condition2p}
\end{eqnarray}
 We are going to use a  Whitney extension. For all $z \in
V$ we denote by $B(z)$ the ball $B(z,10^{-2}d(z,\partial V))$, and
let $X \subset V$ be a maximal set such that for all $z\in X$, the
$B(z)$ are disjoint. Note that by maximality, if $y\in V$, then
$B(y)$ meets some  $B(z)$ for a certain $z \in X$  hence $y \in
4B(z)$ thus the $4B(z)$ cover $V$.

For all $z\in X$ we choose a function $\varphi_z$ which support is
included in $5B(z)$ such that $\varphi_z(y)=1$ for all $y \in
4B(z)$, $0\leq \varphi_z(y)\leq 1$ and $|\nabla \varphi_z(y)|\leq
Cd(z,\partial V)^{-1}$ everywhere. Set $\Phi(y)=\sum_{z\in
X}\varphi_z(y)$ on $V$. We have $\Phi(y)\geq 1$ because  the
$4B(z)$ cover $V$ and the sum is locally finite (because all the
$B(z)$ are disjoint and because the $5B(z)$ that contain a fixed
point $y$ have a radius equivalent to $d(y,\partial V)$. Then we
set $\psi_z(y)=\varphi_z(y)/\Phi(y)$ such that $\sum_{z \in X}
\psi_z(y)=1$ on $V$. Finally, if $m_z$ is the mean value of $v$ on
$B(z)$ we set for all $y \in V$
$$w(y)=\sum_{z\in X}m_z \psi_z(y).$$

If $y \in A_k(1/10)$, $m_z=m_k$ for all $z\in X$ such that $y \in
B(z)$ thus \eqref{condition1p} is verified. In addition,
$$\nabla w (y)= \sum_{z \in X}m_z \nabla \psi_z(z)=\sum_{z \in X}[m_z-m(y)][\nabla \psi_z(y)]$$

where $m(y)$ is the mean value of $v$ on
$B(y)=B(y,10^{-2}d(y,\partial V))$. The sum at the point $y$ has
at most $C$ terms, and all of these terms is less than
$$Cd(y,\partial V)^{-1}|m_z-m(y)|\leq C d(y,\partial
V)^{-3}\int_{10B(y)}|\nabla v|$$ with using Poincaré inequality
and because all the  $5B(z)$ that contain $y$  are included in
 $10B(y)\subset V$. Thus $|\nabla w(y)|\leq C d(y,\partial
V)^{-3}\int_{10B(y)}|\nabla v|$, and to obtain \eqref{condition2p}
it suffice to integrate on $V$, apply Fubini and use
\eqref{condition2}.

Then we apply the co-area formula  (see \cite{f} p.248, and also
\cite{d} chapter 28) to the function $w$ on $V$. We obtain

$$\int_{\R}H^{2}(\Gamma_t)dt=\int_{V}|\nabla w| \leq C \int_V|\nabla u|$$
where $\Gamma_t:=\{y \in V; w(y)=t\}$ is the set of level  $t$ of
the function $w$.  Recall that
$$J(x,r):=r^{-\frac{1}{2}}\min\{\delta_{k,l}; k \not = l\}$$
and
$$\delta_{k,l}=|m_k-m_l|$$
where $m_k$ is the mean value of  $u$ on $D_k$. For all $k_0 \not
= k_1$ we know that $\delta_{k_0,k_1}\geq \sqrt{r}J(x,r)$. Using
Tchebychev inequality we can choose $t_1 \in \R$ such that $t_1$
lies in $\frac{1}{10}[m_{k_0},m_{k_1}]$ and such that
\begin{eqnarray}
H^{2}(\Gamma_{t_1})&\leq& C|m_{k_0}-m_{k_1}|^{-1}\int_{V}|\nabla
u|
\notag \\
&\leq &Cr^{-\frac{1}{2}}J(x,r)^{-1}\int_{V}|\nabla u| \notag \\
&\leq &Cr^{2}J(x,r)^{-1}\omega_2(x,r)^{\frac{1}{2}} \label{choixt}
\end{eqnarray}

For every pair $k_0\not = k_1$ we do the same and choose $t_2$
etc, as many as required by the number of connected components of
$B(x,r)\backslash Z$ (one if $Z$ is a plane, two if $Z$ is a $\Y$
and three if $Z$ is a $\T$). Then we set
$$F=\bigcup_{i}\Gamma_{t_i}\cup [K \cap B(x,r)]\subset B(x,r).$$
$F$ is a closed set in $B(x,r)$ because each $\Gamma_{t_i}$ is
closed in $V=B(x,r)\backslash K$ and $K$ is also a closed set.
Since we have choosing some level sets, $F$ separates the
$A_k(1/10)$ to each other in $B(x,r)$. Indeed, if it is not the
case then there is $k$, $l$ and a continuous path $\gamma$ that
join $A_k(1/10)$ to $A_l(1/0)$ and that is not meeting $K$
(because $K\subset F$). Then $\gamma \subset V$, thus $w$ is well
defined and continuous on $\gamma$, it follows that there is a
point $y \in \gamma $ such that $w(y)={t_i}$. Then, $y \in F$, and
this is a contradiction.

Now we have to prove the $\star$ property. Let  $B(\bar y,s)$ be a
ball centered on $K$ such that $\beta(\bar y,2^ls)\leq
\varepsilon_0 $ for all $0\leq l \leq L$ while  $L$ is the first
integer such that
 $B(\bar y,2^Ls)$ is not included in $B(x,r)$. Set
$B_l:=B(\bar y,2^Ls)$ and possibly by extracting a subsequence we
may suppose using Lemma \ref{recentrage} that in each $B_l$ the
minimal cone associated is almost centered. The radius of $B_l$ is
not as before exactly $2^ls$ but is equivalent with a factor $4$.
Thus the balls $B_l$ forms a sequence of balls centered at
 $\bar y$ such that $B_{l}\subset B_{l+1}$ and $B_0=B(\bar y,s)$. Denote by $Z_l$
the cone associated to $B_l$. We want to show that $F\cap
B(y,s)\subset Z_0(\varepsilon_0):=\{z; d(z,Z_0)\leq \varepsilon_0
s\}$. By definition of $F$, it suffice to show that for all $i$
\begin{eqnarray}
w(y) \not = t_i  \text{ in } B(\bar y,s)\backslash
Z_0(\varepsilon_0). \label{amontrer}
\end{eqnarray}

So let $y \in B(\bar y,s) \backslash  Z_0(\varepsilon_0)$ and
recall that
$$w(y)=\sum_{z \in X}m_z \varphi_z(y).$$
Let $X(y)\subset X$ be the finite set of $z$ such that
$\varphi_z(y)\not = 0$. We claim that
\begin{eqnarray}
\forall z \in X(y), \quad \quad |m_z-m_{D_k}|\leq
Cr^{\frac{1}{2}}\omega_2(x,r)^{\frac{1}{8}} \label{aprouver}
\end{eqnarray}
\begin{sloppypar}
where $m_{D_k}$ is the mean value of $u$ in the appropriate domain
$D_k$ and $m_z$ is the mean value of  $v$ on
$B_z:=B(z,10^{-2}d(z,\partial V))$. First of all, we can use the
proof of Lemma 15 in \cite{l3} to  associate to each connected
component of $B_l\backslash Z_l(\varepsilon_0)$, a component of
$B_{l+1}\cap \{y; d(y,Z_{l+1})\geq 10\varepsilon_0r_l\}$, and by
this way we can rely each component of $B_l\backslash
Z_l(\varepsilon_0)$ to a certain $A_k$ (that contain a $D_k$) (the
argument is just to do an iteration on the scale since we know
that the set $K$ is close the a minimal cone at each scale that we
look at). We denote by $O_0$ the component of $B_s \cap \{y;
d(y,Z_{0})\geq \varepsilon_0s\}$ that contains $y$ and by
induction we denote by $O_l$ the component of $B_l\backslash
Z_l(\varepsilon)$ that is relied to $O_0$. With help of the
particulary geometrical configuration in each $B_l$ we can choose
a domain $G_l$ included at the same time in $O_l$ and in
$O_{l+1}$, and of diameter equivalent to the diameter of $B_l$. We
denote by $m_l$ the mean value of $v$ on $G_l$. We are now ready
to estimate
\end{sloppypar}
\begin{eqnarray}
|m_0-m_L|&\leq& \sum_{l=0}^{L}|m_l-m_{l+1}| \leq
\sum_{l=0}^{L}\frac{1}{|O_{l}|}\int_{O_{l}}|v-m_{l+1}| \notag\\
&\leq&
\sum_{l=0}^{L}C\frac{1}{(2^ls)^{3}}\int_{O_{l+1}}|v-m_{l+1}|
\leq\sum_{l=0}^{L}C(2^{l}s)^{-2}\int_{O_{l+1}}|\nabla v| \notag \\
&\leq&\sum_{l=0}^{L}C(2^{l}s)^{-\frac{1}{2}}
\left(\int_{O_{l+1}}|\nabla v|^2\right)^{\frac{1}{2}} \notag \\
&\leq& \sum_{l=0}^{L}C(2^{l}s)^{-\frac{1}{2}}
\left(\int_{O_{l+1}}|\nabla
v|^2\right)^{\frac{3}{8}}\left(\int_{O_{l+1}}|\nabla
v|^2\right)^{\frac{1}{8}} \notag \\
 &\leq& \sum_{l=0}^{L}C(2^{l}s)^{+\frac{1}{4}}
\left(\int_{O_{l+1}}|\nabla v|^2\right)^{\frac{1}{8}}\label{genial}\\
 &\leq & C\left(\int_{V}|\nabla
v|^2\right)^{\frac{1}{8}}\sum_{l=0}^{L}(2^{l}s)^{\frac{1}{4}}
\leq C\left(\int_{V}|\nabla v|^2\right)^{\frac{1}{8}} \sum_{l=0}^{L} (2^{-l}r)^{\frac{1}{4}} \notag \\
&\leq& C \left(\int_{V}|\nabla v|^2\right)^{\frac{1}{8}}
\sum_{l=0}^{+\infty} (2^{-l}r)^{\frac{1}{4}} \leq Cr^{\frac{1}{4}}
(\int_{V}|\nabla v|^2)^{\frac{1}{8}} \leq
C r^{\frac{1}{4}} (\int_{V}|\nabla u|^2)^{\frac{1}{8}} \notag \\
&\leq & Cr^{\frac{1}{2}}\omega_2(x,r)^{\frac{1}{8}} \label{hl}
\end{eqnarray}

for \eqref{genial} we used the classical estimate on the gradient
of a Mumford-Shah minimizer that is
\begin{eqnarray}
\int_{B(0,R)\backslash K}|\nabla u|^2dx \leq C_N(1+h(R))R^{N-1}
\label{campa}
\end{eqnarray}
obtained by comparing $(u,K)$ and $(v,K')$ where $v$ is equal to
$0$ in $B(0,R)$ and $K'=(K\backslash B(0,R))\cup \partial B(0,R)$.
With the same proof of \eqref{hl} we get
$$|m_L-m_{D_k}|\leq Cr^{\frac{1}{2}}\omega_2(x,r)^{\frac{1}{8}}.$$

On the other hand,  since $z\in X(y)$, then $\varphi_z(y)$ is not
equal to zero. This implies that $d(z,
\partial V)\geq 2d(y,
\partial V)\geq 2 \varepsilon_0s$ thus $B_z:=B(z,10^{-2}d(z,\partial V))\subset
Z_0(\varepsilon_0)^c$. Since by hypothesis  $K$ does not meet this
region,  we can apply Poincaré inequality to prove that
$$|m_z-m_0|\leq Cr^{\frac{1}{2}}\omega_2(x,r)^{\frac{1}{8}}.$$
Finally
$$|m_z-m_{D_k}|\leq |m_z-m_0|+|m_0-m_L|+|m_L-m_{D_k}|\leq Cr^{\frac{1}{2}}\omega_2(x,r)^{\frac{1}{8}}$$
and this completes the proof of \eqref{aprouver}.

Now since $\sum_z\varphi_z(y)=1$ we deduce that

\begin{eqnarray}
|w(y)-m_k|=|w(y)-\sum_{z\in X(y)}\varphi_z(y) m_k|\leq \sum_{z\in
X(y)}|m_z-m_{D_k}|\leq
Cr^{\frac{1}{2}}\omega_2(x,r)^{\frac{1}{8}}.\label{valab}
\end{eqnarray}

Now if we return to the choice of the $t_i$ (see near
\eqref{choixt}) we have taken $t_i \in
\frac{1}{10}[m_{k_0},m_{k_1}]$ for some $k_0$ and $k_1$.

 So thanks to \eqref{valab}, if
$\omega_2(x,r)^{\frac{1}{8}}$ is small enough with respect to
$J(x,r)$ then we are sure that $w(y)\not = t_i$ thus $F$ does not
meet the region $Z_s(\varepsilon_0)$.

We have now to prove \eqref{pr1}. With use of \eqref{pr4} and
\eqref{choixt} we can find a cover of $F$ by a family of balls
$B\{(x_j,r_j)\}$ centered on $K$ and such that $r_j=
C\sqrt{\varepsilon} r$, otherwise we would have a hole in $K$ of
size greater than $C\varepsilon r^2$ which is in contradiction
with \eqref{choixt}. On the other hand, since $\beta_K(x,r)\leq
\varepsilon $, we have $\beta_K(x_j,r_j)\leq C\sqrt{\varepsilon}$.
Now, for every $y\in F\cap B(x_j,r_j)$ we have
$$d(y, Z)\leq d(y,x_j)+d(x_j,Z)\leq C\sqrt{\varepsilon}r+\varepsilon r \leq C \sqrt{\varepsilon}r$$
and the conclusion follows.\qed


Lemma 7 on page 301 of \cite{d} shows how the normalized jump
decreases. So we want to generalize this result to the cones of
type $\T$ and $\Y$ as well. There is no difficulty to do that. We
just have to be careful with the generalized definition of the
jump that depends on the existence of almost centered cones, but
this is not so troublesome. So if the lector already knows how to
control the jump in dimension 2, and if he is convinced that it is
also true for cones of type $\Y$ and $\T$ in $\R^3$, he could just
skip the proofs of the two following lemmas.

\begin{lem} \label{lemme1}Let $(u,K)$ be a Mumford-Shah minimizer in $\Omega$.
Let $x \in K$, $r$ and $r_1$ being such that $B(x,r) \subset
\Omega$ and $0<r_1\leq r \leq \frac{4}{3}r_1$. Suppose in addition
that $\beta(x,r)\leq 100^{-1}$. Then
\begin{eqnarray}
\left|\left(\frac{r_1}{r}\right)^{\frac{1}{2}}J(x,r_1)-J(x,r)\right|\leq
 C\omega_2(x,r)^{\frac{1}{2}}\leq
C(1+h(r)) \label{trot}
\end{eqnarray}
with a constant $C$ that depends only on $N$.
\end{lem}
{\bf Proof :} For all $r_1\leq t\leq 2r_1$ we denote by $Z_t$ a
minimal cone such that
$$\forall y \in K\cap B(x,t), d(y,Z_t)\leq t\beta(x,t)$$
and for all $\lambda$  we also set
$$A_t(\lambda):=\{y \in B(x,t), d(y,Z_t)\geq \lambda t\}.$$
Finally we denote by  $A^k_t$ for $k \in \N\cap [1,\mathpzc{k}_t]$
the different connected components of $A_t$.

To begin, suppose that $Z_{r}$ is almost centered. Recall that in
this case
$$J(x,r)=r^{-\frac{1}{2}}\min\{\delta_{k,l}\}$$
where $\delta_{k,l}=|m_k(r)-m_l(r)|$ and $m_k$ is the mean value
of $u$ on a domain $D_k(x,r)$ in $A_{r}^k(\frac{1}{100})$. Since
$r_1\geq \frac{3}{4}r$ and since $Z_{r}$ is almost centered and
that $\beta(x,r)\leq \frac{1}{100}$, we may consider some balls
$\tilde{D}_k$ in each connected components of $B(x,r_1)\backslash
Z_{r}(\frac{1}{100})$ such that the radii of $\tilde{D}_k$ are
equivalent to $r$ (and thus equivalent to $r_1$) and such that the
$\tilde{D}_k$ are included in $A_r^k(\frac{1}{100})$. By Poincaré
inequality we have
$$|m_{\tilde{D}_k}-m_{A_r^k}|\leq C r^{2} \int_{A_r^k}|\nabla u|$$
and also
$$|m_{D_k(x,r)}-m_{A_r^k}|\leq C r^{2} \int_{A_r^k}|\nabla u|$$
where $m_{D_k(x,r)}$, $m_{A_r^k}$, $m_{\tilde{D}_k(x,r)}$ are the
mean values of $u$ on $D_k(x,r)$, ${A_r^k}$, $\tilde{D}_k(x,r)$.
We deduce that
$$|m_{\tilde{D}_k}-m_{D_{k}(x,r)}|\leq C r^{2} \int_{B(x,r)\backslash K}|\nabla u|\leq Cr^{\frac{1}{2}}\omega_{2}(x,r)
\leq C(1+h(r))r^{\frac{1}{2}}.$$ The last inequality comes from
\eqref{campa}. By the same way we obtain
$$|m_{\tilde{D}_k}-m_{D_{k}(x,r_1)}|\leq C r^{2}
\int_{B(x,r)\backslash K}|\nabla u|\leq
Cr^{\frac{1}{2}}\omega_{2}(x,r)\leq C(1+h(r))r^{\frac{1}{2}}$$
where the $D_{k}(x,r_1)$ are the domains in the definition of
$J(x,t_1)$. This then gives the estimation of
$r^{\frac{1}{2}}J(x,r)-r_1^{\frac{1}{2}}J(x,r_1)$ to prove
\eqref{trot}.

Finally if $Z_r$ is not almost centered then we have two cases.
The first one is when $Z_{r_1}$ is nether almost centered and then
we can use  $2r_1$ and $2r$ and that is the same as the above
argument. The second case is when $Z_{r_1}$ is almost centered and
then this implies that $Z_{r_1}$ is a cone of minor type than the
type of $Z_r$ thus it suffice to control the mean values only in
connected components $A_r$ that meets the $A_{r_1}$, and the
difference between those mean values are always bounded by the
jump $J(x,r)$.\qed

\begin{lem} \label{lemme7} Let $(u,K)$ be a Mumford-Shah minimizer in $\Omega$.
Then if $x \in K$ and $r$ are such that $B(x,r) \subset \Omega$
and for all $r_1<t<r$, $\beta(x,t)\leq 10^{-1}$, then
\begin{eqnarray}
J(x,r_1)\geq
\left(\frac{r}{r_1}\right)^{\frac{1}{2}}\left[J(x,r)-C'\right]
\label{ouf}
\end{eqnarray}
where $C':=C(1+h(r))$ and $C$ depends only on $N$.
\end{lem}

{\bf Proof :} If $r_1\leq r \leq \frac{4}{3}r_1$ then \eqref{ouf}
is a consequence of Lemma \ref{lemme1}. Otherwise we use a
sequence of radii $r_k$ such that $r_k=\frac{4}{3}r_{k-1}$ and we
apply Lemma \ref{lemme1} a number of time until  $r_k$ is greater
than $r$. We obtain
\begin{eqnarray}
J(x,r_1)&\geq& \sqrt{4/3}^{k}J(x, (4/3)^{k}r_1)-C
\sqrt{4/3}^k(1+\sqrt{4/3}^{-1}+\sqrt{4/3}^{-2}+\dots)
\notag \\
&\geq
&\sqrt{4/3}^k[J(x,(4/3)^kr_1)-\frac{C\sqrt{4/3}}{1-\sqrt{4/3}}].
\end{eqnarray}
and we conclude by using Lemma \ref{lemme1} a last time.\qed

In the following we will sometimes use the notations $F$ and $B$
instead of $F(x_0,r_0)$ and $B(x_0,r_0)$. In addition, without
loss of generality we may suppose by now that  $x_0=0$.


\subsection{Stopping times balls and bad mass}

Our goal in this section is to construct a family of balls   $S$
by a stopping time argument, with the condition that in all balls
of $S$,  the singular set $K$ will always looks like a minimal
cone.

We suppose that $B(0,4r_0)\subset \Omega$. For all $x \in K \cap
B(0,r_0)$ and $r>0$, we say that $B(x,r)$ is a good ball (and then
 denote $B(x,r)\in \mathcal{G}$) if
\begin{eqnarray}
H^2(F\cap B(x,r))-H^2(K\cap  B(x,r)) \leq \varepsilon_0' r^2
\label{bon1}
\end{eqnarray}
and also if there is a minimal cone $Z$ such that
\begin{eqnarray} \forall y \in K\cap  B(x,r),\;
d(y,Z)\leq \varepsilon_0 r .\label{bon2}
\end{eqnarray}

Here, $\varepsilon_0'$ and $\varepsilon_0$ are such that
$\varepsilon < \varepsilon_0' < \varepsilon_0<10^{-5}$. Note that
since $\beta(0,4r_0)\leq \varepsilon$ the radii of balls that
don't verify \eqref{bon2} is bounded by
$\frac{\varepsilon}{\varepsilon_0}r_0$ and if
$C\omega_2(0,r_0)^{\frac{1}{2}}J(0,r_0)^{-1}\leq \varepsilon$, the
radii of balls that don't verify \eqref{bon1} is bounded by
$\sqrt{\varepsilon}r_0$.

Now, for all $x \in K$ we define the stopping time function
$$d(x):= \inf\{r; \forall t\geq r, B(x,t)\in \mathcal{G}\}.$$

Then with help of the Vitali covering lemma, from the collection
of balls $$\{B(x,Ad(x))\}_{x \in K\cap B(0,r_0)}$$ with $A$ a
constant that will be chosen later, we get a disjoint subfamily
$\{B_i\}_{i \in I}$ such that $\{5B_i\}_{i \in I}$ is covering.
Denote $S:=\{B_i\}_{i \in I}$  the ``Bad balls''.
 For all $r\leq r_0$ we set
$$I_r:=\{i \in I ; B_i\cap B(0,r) \not = \emptyset\}$$
and we introduce a new quantity called ``Bad mass'' defined by
$$m(0,r):=\frac{1}{r^2}\sum_{i \in I_r}r_i^2.$$

By convention, a single point $\{x\}$ with $d(x)=0$ will be
identified with the ball $B(x,d(x))$.

\subsection{Whitney extension}

Here we have to recall some definitions and a result from
\cite{l3} in a little weaker form.

Let $K$ be a closed set in $\bar B(x_0,r_0)$ such that
$H^{2}(K\cap \bar B(x_0,r_0))< + \infty$. Suppose that there is a
positive constant $\varepsilon_0 <10^{-5}$ and a minimal cone $Z$,
centered at $x_0$, such that
\begin{eqnarray}
\sup \{ d(x,Z) ; x \in K\cap B(x_0,r_0)\} \leq r_0\varepsilon_0
\label{sect1}
\end{eqnarray}
 and that $K$ is
separating in $B(x_0,r_0)$. For all $x \in K\cap B(x_0,r_0)$ and
$r>0$ such that $B(x,r)\subset B(0,r_0)$ recall that
$$\beta(x,r)=\inf_{Z \ni x}\frac{1}{r}\sup \{d(x,Z); x \in K\cap B(x,r)\}.$$
Let $\rho \in [\frac{1}{2}r_0,\frac{3}{4}r_0]$ and assume that we
have an application
\begin{eqnarray}
\delta: B(x_0,\rho) \to [0,\frac{1}{4}r_0] \label{sect12}
\end{eqnarray}
 with the property that
\begin{eqnarray}
\beta(x,r)\leq \varepsilon_0, \text{ for all } x \in K\cap
B(x_0,\rho) \text{ and } r \text{ such that  } \delta(x)\leq r
\leq \frac{1}{4}r_0 \label{sect13}.
\end{eqnarray}
In addition we suppose that
\begin{eqnarray}
\delta \text{ is } C_0-\text{Lipschitz}\label{sect14}.
\end{eqnarray}
 The application $\delta$ will be called the ``geometric
function''.

\begin{defin}[Hypothesis $\mathcal{H}$] \label{defhyph} We will say that a closed set $K\subset B(x_0,r_0)$ with finite
$H^2$ measure is satisfying hypothesis
$\mathcal{H}$ if \\
i) There is a minimal cone $Z$ that verify \eqref{sect1} for a
``geometric constant"
 $\varepsilon_0<10^{-5}$ and a ``Lipschitz constant" $C_0$. \\
ii)  $K$ is separating in $B(x_0,r_0)$.\\
iii) There is a geometric function $\delta$ satisfying
\eqref{sect12}, \eqref{sect13} and \eqref{sect14} for a radius
$\rho \in [\frac{1}{2}r_0, \frac{3}{4}r_0]$.
\end{defin}


 Let $U>1$ be a constant that will be fixed later, depending on $C_0$
and a dimensional constant. In addition we assume that
$\varepsilon_0$ is very small compared to $U^{-1}$. For all $t>0$
we define
\begin{eqnarray}
\mathpzc{V}:= \bigcup_{x \in K\cap B(0,\rho)}B(x,
\frac{10}{U}\delta(x)). \label{defmore}
\end{eqnarray}
We also set
\begin{eqnarray}
\mathpzc{V}_{\rho}:= \bigcup_{x ; B(x, \frac{10}{U}\delta(x))\cap
\partial B(x_0,\rho)\not = \emptyset}B(x, \frac{10}{U}\delta(x)). \label{defQ3}
\end{eqnarray}

Recall that  by hypothesis, $K$ is separating in $B(x_0,r_0)$ and
that for all $k \in [1,\mathpzc{k}^{B(x_0,r_0)}]$ we denote by
$A_k(x_0,r_0)$ the connected components of $B(x_0,r_0)\backslash
Z_{\varepsilon_0r_0}$ and by $\Omega_k(x_0,r_0)$ the connected
component of $B(x_0,r_0)\backslash K$ that contains
$A_k(x_0,r_0)$. We also set
\begin{eqnarray}
\Delta_k:=B(x_0,\rho)\cap (\Omega_k(x_0,r_0)\cup \mathpzc{V}).
\label{defQ4}
\end{eqnarray}

Then we have the following lemma.

\begin{lem}{\rm \cite{l3}}\label{extentionW}{\rm(Whitney Extension)} Let $K$ be a closed
set in $B(x_0,r_0)$ satisfying Hypothesis $\mathcal{H}$ with a
geometric function $\delta$, a minimal cone $Z$, a constant
$\varepsilon_0<10^{-5}$ and a radius $\rho\in [\frac{1}{2}r_0,
\frac{3}{4}r_0]$. Then for all function $u \in
W^{1,2}(B(0,r_0)\backslash K)$, and for all $k \in
[1,\mathpzc{k}^{B(x_0,r_0)}]$, there is a function
$$v_k \in
W^{1,2}(\Delta_k \backslash \mathpzc{V}_{\rho})$$ such that
$$v_k = u \text{ in }B(x_0,\rho)\backslash \mathpzc{V}$$
and
\begin{eqnarray}
\int_{\Delta_k\backslash \mathpzc{V}_{\rho}}|\nabla v_k|^2dx \leq
+ C\int_{B(0,r_0)\backslash K}|\nabla u|^2dx \label{inwh1}
\end{eqnarray}
where $C$ is a constant depending only on dimension and where
$\mathpzc{V}$, $\mathpzc{V}_{\rho}$, and $\Delta_k$ are defined in
\eqref{defmore}, \eqref{defQ3}, and $\eqref{defQ4}$ with constant
$U>30C_0$ depending also on dimension.
 \end{lem}

From balls of $S$,  we want to apply Lemma \ref{extentionW} to get
a good extension of $u$ near the bad balls. This extension will
allow us replace in each bad ball the set $K$ by a new set in
order to get some estimates. So we begin by introducing a
geometric function associated to the balls of $S$. We define
\begin{eqnarray}
\forall x \in \R^3; \quad \delta(x):=\inf_{B\in  S}\{d(x, B)+r_B\}
\label{defd}
\end{eqnarray}
where $r_B$ is the radius of the ball $B$ (that could be equal to
$0$).

\begin{prop}\label{hypothesisH}
Application $\delta$ is a geometric function associated to $F$ in
$B(0,r_0)$ for all $\rho\in [\frac{1}{2}r_0,\frac{3}{4}r_0]$, with
Lipschitz constant $C_0=1$ and geometric constant
$10\varepsilon_0$. In addition, we have Hypothesis $\mathcal{H}$
on $F$ in $B(x_0,r_0)$ and
$$\bigcup_{i \in I}\frac{10}{U}B_i \subset \mathpzc{V}$$
where $\mathpzc{V}$ is defined in \eqref{defmore}.
\end{prop}

\begin{rem}  Note that
since $C_0=1$, $U$ is depending only on dimension.
\end{rem}

{\bf Proof :} We have to verify \eqref{sect12}, \eqref{sect13} and
\eqref{sect14}. Let $\rho \in [\frac{1}{2}r_0,\frac{3}{4}r_0]$.
Recall that $F$ is separating in $B(0,r_0)$ and
$$F\cap B(0,r_0)\subset \{y \in B(x_0,r_0); d(y,Z^0) \leq r_0C\sqrt{\varepsilon}\}.$$
Then if $\varepsilon$ is small enough with respect to
$\varepsilon_0$, for all $x \in F\cap B(0,\rho)$ and for all ball
$B(y,r)$ with $r>\frac{1}{100}r_0$ that contain $x$ we have
$\beta(y,r)\leq \varepsilon_0$ thus for all $x \in F \cap
B(0,\rho)$ we easily have
$$\delta(x)\leq \frac{1}{4}r_0$$
and \eqref{sect12} is proved.

 Now let $x\in F \cap B(0,\rho)$ and let $r$ be a radius such that $\delta(x)\leq r \leq
 \frac{1}{4}r_0$.  Let $B$ be a ball of $S$ such that
 \begin{eqnarray}
d(x,B)+r_B \leq 2r \notag
\end{eqnarray}
(we always have one $B$ like that by definition of $\delta$). Let
$x_B$ be the center of $B$. Then we have
$$d(x,x_B)\leq  2r$$
thus $x \in B(x_B,2r)$ and $B(x,r) \subset B(x_B,3r)$. Since $r
\geq r_B \geq d(x)$, we know by definition of $d(x)$ that
$\beta(x_B,3r)\leq \varepsilon_0$. Moreover, for all $t >3r$ we
have
$$\beta(x_B,t)\leq \varepsilon_0.$$
Then we can apply Property $\star$ in $B(x_B,3r)$ in order to get
a cone $Z$ containing $x_B$ such that for all $y \in F\cap B(x_B,
3r)$, $d(y,Z)\leq  \varepsilon_0 3r$. Since $B(x,r) \subset
B(x_B,3r)$ and $x \in F$ we deduce that $\beta_F(x,r)\leq
10\varepsilon_0$ and \eqref{sect13} is proved with
$10\varepsilon_0$ instead of $\varepsilon_0$.

Finally, if $B$ is a ball of $S$ then for all $x$ and $y$ we have
\begin{eqnarray}
d(x,B)&\leq& d(x,y)+d(y,B) \notag \\
d(x,B)+r_B&\leq& d(x,y)+ d(y,B)+r_B \notag \\
\delta(x) \leq d(x,B)+r_B &\leq& d(x,y)+d(y,B)+r_B \notag \\
\delta(x) &\leq & d(x,y)+ d(y,B)+r_B \notag
\end{eqnarray}
then passing to the infimum we deduce
$$|\delta(x)-\delta(y)|\leq d(x,y)$$
and that shows that the application $x\mapsto \delta (x)$ is
$1$-Lipschitz.

So we deduce that we have Hypothesis $\mathcal{H}$ on $F$ in
$B(x_0,r_0)$ with application $\delta$ defined in \eqref{defd}.
Let us show that
\begin{eqnarray}
\bigcup_{i \in I}\frac{10}{U}B_i\subset \mathpzc{V} \label{boucle}
\end{eqnarray}

Let $B_i=B(x_i,r_i)\in S$ be a bad ball. We claim that
$$\delta(x_i)\geq r_i.$$
Indeed, recall that the balls $B \in S$ are disjoint. If we take
$B_i$ in the infimum of the definition of $\delta$ we get $d(x_i,
B_i)+r_i=r_i$ and if we take a ball $\tilde B$ out of $B$ we get
again $d(x_i,\tilde{B})+r_{\tilde B}\geq r_i$. Thus
\begin{eqnarray}
B(x_i,\frac{10}{U}r_i)\subset B(x_i,\frac{10}{U}\delta(x_i))
\subset \mathpzc{V}. \notag  \qed
\end{eqnarray}

\section{Useful estimates}

We are now ready to compute some estimates about different
quantities that will lead to regularity. The main point  is to
show some decay estimates on the normalized energy
$\omega_2(x,r)$. This decay will come from the same sort of
argument as in \cite{d1}. In dimension 2, the intersection between
$\partial B(x,r)$ and $K$ is mainly constituted of single points.
Here in dimension 3, $\partial B(x,r)\cap K$ is more complicated
and this will led some problems. We start by finding a judicious
radius  $\rho$ to begin the estimates.

\subsection{Choice of the radius}

For the choice of the radius we select a  $\rho  \in
R:=[\frac{r_0}{2}, \frac{3}{4}r_0]$ such that the mass of the bad
balls $\{B_i\}_{i \in I}$ that are meeting  $\partial B(0,\rho)$
is less than average. Recall that the $B_i$ are the bad balls
$B(x_i,r_i)\in S$. Set $I(\rho):=\{i \in I ; B_i \cap
\partial B(0,\rho)\not = \emptyset\}$
   and let $r_i$ be the radius of $B_i$. By such a choice of $\rho$ we have
$$\sum_{i\in I(\rho)}r_i^2 \leq \frac{1}{|R|}\int_{R}\sum_{i \in I(t)} r_i^2 dt\leq \frac{1}{|R|}\sum_{i \in I}
\int_{t ; i \in I(t)}r_i^2\leq C \frac{1}{|R|} \sum_{i \in
I}r_i^3.$$ Finally we have found a $\rho$ that verify
\begin{eqnarray}
\sum_{i\in I(\rho)}r_i^2 \leq \frac{C}{r_0} \sum_{i \in I} r_i^3
\leq C \sup_{i}\{r_i\}\sum_{i \in I}r_i^2 \leq
C\sqrt{\varepsilon}r^2m(r). \label{energ1}
\end{eqnarray}


\subsection{Comparaison with an energy minimizing function}

Since $\rho$ is now chosen, we are ready to compare with an energy
minimizing function and use the decay result of \cite{l3}. By
construction of $S$, the set $F$ is
$(\varepsilon_0,\sqrt{\varepsilon})$-minimal in sense of
Definition 8 of \cite{l3}. In fact, we know that
 $F$ is $\varepsilon_0$-minimal in the complement of the $\{B_i\}_{i\in I}$,
 and for all $i$, we have that $r_i\leq
\sqrt{\varepsilon}r_0$. Set
$$
G:=F^\rho=(F\backslash \bigcup_{i \in I(\rho)} B_i) \cup
\bigcup_{i \in I(\rho)}\partial B_i.
$$

Then if  $\varepsilon$ is small enough with respect to
$\varepsilon_0$ and $\varepsilon_2$ (the constant of \cite{l3}) we
can apply Theorem 9 of \cite{l3}. Thus we know that the normalized
energy decreases for all energy minimizer in $B(0,r_0)\backslash G
$. In particular if $w$ is the energy minimizer in $B\backslash G$
that is equal to $u$ on $\partial B \backslash G=\partial B
\backslash F$ (for the existence of such a minimizer, one can see
for example \cite{d} page 97), applying Theorem 9 of \cite{l3}
with $0<\gamma<0,8$, we have that for all $a<\frac{1}{2}$, there
is a $\varepsilon_2$ (that depends on $a$ and $\varepsilon_0$),
such that if $\varepsilon$ is small enough (depending on
$\varepsilon_0$ and $a$),
\begin{eqnarray}
\frac{1}{(ar_0)^2}\int_{B(0,ar_0)\backslash G}|\nabla w|^2\leq
a^{\gamma}  \frac{1}{r_0^2}\int_{B(0,r_0)\backslash G}|\nabla
w|^2. \label{decr1}
\end{eqnarray}

The second useful fact is the following. Since $(u,K)$ is a
Mumford-Shah minimizer and $(w,G)$ is a competitor we have
\begin{eqnarray}
\int_{B(0,\rho)\backslash K} |\nabla u|^2 + H^2(K\cap B(0,\rho))&
\leq &\int_{B(0,\rho)\backslash G}|\nabla w|^2+ H^2(G\cap
B(0,\rho))+\rho^2h(\rho). \notag
\end{eqnarray}
Hence
\begin{eqnarray}
\int_{B(0,\rho)\backslash K} |\nabla u|^2&-&\int_{B(0,\rho)\backslash G}|\nabla w|^2 \notag \\
&\leq& H^2(G\cap B(0,\rho))-H^2(K\cap B(0,\rho))+\rho^2h(\rho) \notag \\
&\leq& Cr_0^2\omega_2(0,r_0)^{\frac{1}{2}}J(x_0,r_0)^{-1}+C\sum_{i \in I(\rho)}r_j^2+\rho^2h(\rho) \notag \\
&\leq & Cr_0^2\omega_2(0,r_0)^{\frac{1}{2}}J(x_0,r_0)^{-1}+C
\sqrt{\varepsilon}r_0^2m(0,r_0)+\rho^2h(\rho).\label{estimhort}
\end{eqnarray}

The third point is that $\nabla  w$ and $\nabla(w-u)$ are
orthogonal in $L^2(B(0,r_0))$. That comes from the fact that $w$
is an energy minimizer in $B(0,r_0)\backslash G$ and $u$ is a
competitor for $w$. Thus
$$\int_{B(0,r_0)\backslash G}|\nabla  u - \nabla w|^2= \int_{B(0,r_0) \backslash G}|\nabla u|^2 -
 \int_{B(0,r_0)\backslash G}|\nabla w|^2.$$
We can now estimate the energy of $u$. Let $0<a<\frac{1}{2}$, then
\begin{eqnarray}
\int_{B(0,ar_0)\backslash G}|\nabla u|^2&\leq&
2\int_{B(0,ar_0)\backslash G} |\nabla w|^2 +
2\int_{B(0,ar_0)\backslash G}
|\nabla w-\nabla u|^2\notag\\
&\leq& 2a^{2+\gamma} \int_{B(0,r_0)\backslash G} |\nabla w|^2 +
2\int_{B(0,r_0)\backslash G}
|\nabla w-\nabla u|^2 \notag \\
&\leq& 2a^{2+\gamma} \int_{B(0,r_0)\backslash G} |\nabla  u |^2 +
2\int_{B(0,r_0)\backslash G} |\nabla  u|^2
-2\int_{B(0,r_0)\backslash G} |\nabla w |^2 . \notag
\end{eqnarray}
Hence,
\begin{eqnarray}
\omega_2(0,ar_0)\leq 2a^{\gamma} \omega_2(0,r_0) +
C\frac{1}{a^2}\omega_2(0,r_0)^{\frac{1}{2}}J(0,r_0)^{-1}+C
\frac{\sqrt{\varepsilon}}{a^2}m(0,r_0)+\frac{1}{a^2}\rho^2h(\rho)
\label{debdec}.
\end{eqnarray}

Inequality \eqref{debdec} is the fundamental estimate that will be
used to control the energy.


\subsection{Compactness lemmas  for almost minimal sets}

The purpose of this section is to show some geometrical results
about almost minimal sets (see definition \ref{m2}). In the future
estimates, we will use an argument which allows us to win
something in each bad ball, in order to prove that there are not
so many. The main lemma says the following. If  $B(x,r)$ is a ball
such that $x \in K$ and $\beta(x,r)\leq \varepsilon_0$ but
$\beta(x,ar)>\varepsilon_0$, then there is a set that does better
than $K$ in $B(x,r)$ in terms of $H^2$-measure.

Recall that for any almost minimal set $E$ in $B(x,r)$, we denote
by $f(r)$ the excess of density
$$f(r)=\theta(x,r)-\lim_{t\to 0} \theta(x,t)$$
with
$$\theta(x,r)=r^{-2}H^{2}(E\cap B(x,r)).$$
The limit at $0$ of $\theta$ exists  because  $E$ is almost
minimal (see 2.3. of \cite{d3}). For $x\in E$ we call $d(x)$ the
density at $x$, that is $d(x)=\lim_{t\to 0} \theta(x,t)$. The
function $d(x)$ can only take a finite number of values, more
precisely $d(x)\in \{0, \pi, \frac{3\pi}{2}, d_+\}$ that are
(excepted $0$) densities of the three minimal cones in $\R^3$.

For an almost minimal set $E$, the function $\theta(x,t)$ is
 non decreasing in $t$ thus the limit when
 $t$ tend to $0$ exists and that allows us to define the function $d(x)$. Unfortunately, if  $E$
 is now the singular set of a Mumford-Shah minimizer, the monotonicity of
  $\theta$ is not known.  So we have some difficulties so define the
analogue of $f(r)$ for a Mumford-Shah minimizer.

In order to use Theorem \ref{gd3}, we want to control $f(r)$. That
will be the role of the following Lemmas. Our goal is to obtain a
statement analogous to Theorem  \ref{gd3} but with only an
hypothesis on $\beta(0,r_0)$ instead of $f(0,r_0)$.

First of all, an application of Proposition 16.24 of \cite{d3} in
$B(x,r10^{-3})$ with $\eta_1=\varepsilon_210^3$, mixed with
Proposition 18.1 of \cite{d3} in $B(x,r10^{-5})$ and
$\eta_1=\varepsilon_710^{-5}$ (where $\varepsilon_7$ and
$\varepsilon_3$ are defined in \cite{d3}) allows us to state the
following lemma. Recall that $D_{x,r}$ is the normalized bilateral
Hausdorff distance.

\begin{lem}{\rm \cite{d3}} \label{compa1} For each choice of $b\in $There is a $\eta_1\geq 0$ such
that if $E$ is an almost minimal set in an open set $U \in \R^3$,
with gauge function $h(r)=C_0r^{b}$, if $x \in E$ and $r>0$ are
such that $B(x,r)\subset U$, if there is $Z$, centered at $x$, of
type $\mathbb{P}$, $\mathbb{Y}$ or $\mathbb{T}$ such that
$$D_{x,r}(E,Z)\leq \eta_1 ,\quad   \quad h(2r)\leq \eta_1 ,\quad \quad \int_{0}^{2r}h(t)\frac{dt}{t}\leq \eta_1$$
and if $E$ is separating in  $B(x,r)$, then there is a point $x
\in E\cap B(x,r10^{-5})$, of the same type of $Z$.
\end{lem}

We say that $x$ has the same type as $Z$ if $d(x)$ is equal to the
density of the cone $Z$.

\begin{rem}\label{compa1r} The hypothesis of separating
are only useful for the case of $\mathbb{T}$. See Propositions
16.24 and 18.1 of \cite{d3} for more details.
\end{rem}

\begin{rem} Lemma \ref{compa1} is not trivial because we can
imagine  that $E$ is very close to a cone of type $\mathbb{T}$ in
$B(x,r)$ but contains only $\Pp$ points and $\mathbb{Y}$-points
(see \cite{d3} Section 19). The lemma says that under separating
conditions and if $h$ and $\beta$ are small enough, this is not
possible.
\end{rem}

Here is now the statement that will be useful for the next
sections. The reader is invited to compare it with Theorem
\ref{gd3}.

\begin{lem}\label{compa2}For each choice of $b \in (0,1]$, $\bar
c>0$ and  $C_0>0$ we can find $\eta_2>0$ and $C\geq 0$ such that
the following holds. Let $E$ be a reduced MS-almost minimal set in
$\Omega \subset \R^3$ with gauge function $h$. Suppose that $0 \in
E$, $r_0>0$ such that $B(0,110r_0)\subset \Omega$ and $h$ is
satisfying
$$h(r)\leq C_0r^b\quad \text{ for } 0<r<220r_0 $$
and
$$h(220r_0)\leq \eta_2 ,\quad \quad
\int_{0}^{220r_0}h(t)\frac{dt}{t}\leq \eta_2.$$
 Assume in addition that
$$D_{0,100r_0}(E,Z)\leq \eta_2$$
where $Z$ is a minimal cone centered at the origin such that
$$H^2(Z\cap B(0,1))\leq d(0).$$
Then for all $x \in E \cap B(0,4r_0)$ and for all $0<r<5 r_0$
there is a minimal cone $Z(x,r)$ such that
$$D_{x,r}(E,Z(x,r))\leq \bar c\left(\frac{r}{r_0}\right)^{\alpha}.$$
\end{lem}

{\bf Proof :} We take  $\eta_2<\varepsilon_1$ (the constant of
Theorem \ref{gd3}). In order to apply Theorem \ref{gd3}, all we
have to prove is that
$$f(0,110r_0)\leq \varepsilon_1.$$
 If $\eta_2$ is smaller than $\eta_1$ we can apply Lemma
\ref{compa1} to $E$ in $B(x,110r_0)$ thus there is a point $z$ in
$B(x,10^{-3}r_0)$ of same type of $Z$. In particular
$d(z)=H^2(Z\cap B(z,1))=\frac{1}{r^2}H^2(Z\cap B(z,r))$ for all
$r$. Hence we can compute the excess of density at $z$ in
$B(z,55r_0)$ by
$$f(z,55r_0)=\frac{1}{(55r_0)^2}[H^2(E\cap B(z,55r_0))-H^2(Z\cap B(z,55r_0))].$$
Now define a competitor $L$ by
$$L= \left\{
\begin{array}{cc}
M \cup Z\cap B(z,55r_0) & \text{ in } \bar B(z,55r_0) \\
E & \text{ in } \Omega\backslash B(z,55r_0)
\end{array}
\right.
$$
where $M$ is a little wall:
$$M:=\{x \in \partial B(z,55r_0) ; d(x,Z)\leq 500\eta_2r_0\}.$$
The set $L$ is a MS-competitor for $E$ thus
\begin{eqnarray}
H^2(E\cap B(z,55r_0))&\leq& H^2(L\cap B(z,55r_0))+(55r_0)^2h(55r_0) \notag \\
&\leq & H^2(M)+H^2(Z\cap B(z,55r_0))+(55r_0)^2h(55r_0). \notag
\end{eqnarray}
Since $H^2(M)\leq Cr_0^2\eta_2$ we deduce
$$f(z, 55r_0)\leq C\eta_2.$$
Now if $\eta_2$ is small enough  compared to $\varepsilon_1$, we
can apply Theorem \ref{gd3} in $B(z,55r_0)$ then for all $y \in
E\cap B(z,5r_0)$ and $0<r<10r_0$ we have
\begin{eqnarray}
\beta(y,r)\leq \bar c\left( \frac{r}{r_0}\right)^\alpha
\label{estimdec}.
\end{eqnarray}
In addition, since $d(x,z)\leq 10^{-3}r_0$ we deduce that
\eqref{estimdec} is true for all $y \in B(x,4r_0)$ and
$0<r<10r_0$.\qed

\begin{defin}\label{remr} By now we will call $\bar \eta_2$ the constant given by Lemma \ref{compa2}
with $\bar c=1$, $r_0=\frac{1}{200}$  $C_0=0$ and $b=0$, and we
call $\bar r$ the radius such that
$$\left(\frac{\bar r}{2}\right)^\alpha =\frac{1}{2}\varepsilon_0.$$
\end{defin}

Now we are ready to prove our fundamental lemma that will be used
later to count the mass of bad balls.

\begin{lem} \label{compacite}
For all $\varepsilon_0>0$, and  for all $r<\bar r$, there is a
constant $\eta_0$ such that if $E$ is a closed set of finite $H^2$
measure in $B(0,1)\subset \R^3$ that contains the origin, with the
uniform concentration Property (with constant $C_u$), and assume
that
\begin{eqnarray}
   \beta(0,1)&\leq &\bar \eta_2 \label{beta1} \\
\beta(0,r)&\geq & \varepsilon_0 \label{beta2}
\end{eqnarray}
such that the cone in $\beta(0,1)$ is centered in $B(0,10^{-5})$.
 If in addition we assume that there is a set $F$ that contains
 $E$, that is separating in $B(0,1)$ (see Definition
\ref{separ}) and such that
$$H^2(F)-H^2(E)\leq \eta_0.$$
 Then there is a MS-competitor $L$ for $E$ in $B(0,\frac{3}{4})$ such that
$$H^2(E)-H^2(L) \geq \eta_0$$
\end{lem}
{\bf Proof : } The argument is by contradiction. If the lemma is
not true, then there is a $r<\bar r$ and there is an
$\varepsilon_0<\frac{1}{100}\bar \eta_2$ such that for all $\eta$
there is a set $E_\eta$ that verify \eqref{beta1} and
\eqref{beta2}. In addition for all MS-competitor $L_\eta$ we have
\begin{eqnarray}
H^2(E_\eta)- H^2(L_\eta)\leq \eta. \label{minsurf}
\end{eqnarray}
 And for all $\eta$ there is a set $F_\eta$ that contain $E_\eta$, is separating in
 $B(0,1)$,
 and such that
\begin{eqnarray}
H^2(F_\eta)-H^2(E_\eta)\leq \eta. \label{separ}
\end{eqnarray}
Now let $\eta$ tend to $0$. Passing if necessary to a subsequence,
we may assume that the sequence
 of sets $E_\eta$ converges to a certain $E_0$ in sense of Hausdorff
 distance. Passing to the limit, we deduce that this set
$E_0$ still verify \eqref{beta1} and \eqref{beta2} .

We want to show that $E_0$ is a minimal set in $B(0,\frac{3}{4})$.
Let $L$ be a MS-competitor for $E_0$ in $B(0,\frac{3}{4})$. Since
$E_\eta$ tend to $E_0$ for the
  Hausdorff distance $D_H$, we know that for all $\tau$ there is a $\eta'$
  such that for all
    $\eta < \eta'$, $D_H(E_0,E_\eta)\leq \tau$. Thus if
   $T_\tau:=\{x \in \partial B(0,1); d(x,L)\leq \tau \}$, we
   deduce that
   $E_\eta
   \cap \partial B(0,1)\subset T_\tau$. Therefore, the set
   $L_\eta:=L\cup (E_\eta\cap B(0,1)\backslash B(0,\frac{3}{4}))\cup T_\tau$
   is a  MS-competitor for $E_\eta$. Then applying  \eqref{minsurf}
   we obtain
\begin{eqnarray}
H^2(E_\eta\cap B(0,\frac{3}{4}))&\leq& H^2(L_\eta)+ \eta \notag \\
& \leq& H^{2}(L)+ H^{2} (T_\tau) +\eta \notag \\
&\leq &  H^{2}(L\cap B(0,\frac{3}{4}))+ \eta+C\tau .\notag
\end{eqnarray}

 In addition, by hypothesis the sets $E_\delta$ verify the uniform concentration
 property with same constant $C_u$. This allows us to say that (see \cite{d} section 35)
\begin{eqnarray}
H^2(E_0)\leq \underline{\lim}_{\eta\to 0} H^2(E_\eta) .\notag
\end{eqnarray}
Hence, letting $\eta$ tend to $0$ we obtain
$$H^{2}(E_0)\leq H^2(L)+\tau$$
then letting $\tau$ tend to $0$,
$$H^{2}(E_0)\leq H^2(L)$$
 thus $E_0$ is a minimal set (i.e.
almost minimal set with gauge function equal to zero).

On the other hand, $E_0$ is separating in  $B(0,1)$, because if it
is not the case, we can find a continuous path $\gamma$ that join
$A^+$ and $A^-$ (two points in different connected component of
$B(0,1)\backslash Z_{10^{-5}}$) in $B(0,1)$ and such that $\gamma$
does not meet  $E_0$. Since $E_\eta$ converge to $E_0$ for the
Hausdorff distance, for all $\tau$ there is a $\eta_\tau$ such
that if $\eta < \eta_\tau$, all the $E_\eta$ are  $\tau$ close to
$E_0$. Let $x$ be the point of $\gamma$ that realize the infimum
of $d(x,E_0)$. Since $\gamma$  is disjoint from $E_0$, there is a
ball centered at $x$ with positive radius  $r$ that is not meeting
$E_0$. Thus if we choose $\eta$ smaller than $r$ we get that all
the $E_\eta$ for $\eta < \eta_\tau$ contain a hole of size $r$,
but this is not possible according to \eqref{separ}.

Thus finally  $E_0$ is a minimal set in $B(0,\frac{3}{4})$, which
is separating and verifies \eqref{beta1} and \eqref{beta2}. We
want now to apply Lemma \ref{compa2} to obtain a contradiction. We
know that
$$\beta(0,1)\leq \bar \eta_2$$
and that the cone associated is centered in $B(0,10^{-5})$.  We
claim that
\begin{eqnarray}
D_{z,\frac{1}{2}}(E_0,Z) \leq \bar \eta_2. \label{deil}
\end{eqnarray}
All we have to show is that for all $x \in Z$, $d(x,E_0) \leq \bar
\eta_2$. If it is not the case, then we can find  $x \in Z$ such
that $B(x,\bar \eta_2)\cap E_0=\emptyset$. But then we can find a
continuous path that join two different connected components of
 $B(0,1)\backslash
Z_{\frac{1}{100}\eta_2}$ without meeting $E$, and that is not
possible if $E$ is separating. So we have shown \eqref{deil} and
then we can apply Lemma \ref{compa2} in $B(z,\frac{1}{2})$ (i.e.
$r_0=\frac{1}{200})$, which implies that
$$\beta(0,r)\leq \frac{1}{2}\varepsilon_0$$
because of the definition of $\bar r$, and this yields a
contradiction with \eqref{beta2} so the proof is now complete.
\qed

Applying Lemma \ref{compacite} we can deduce to following
proposition.

\begin{prop}\label{compa} Let $i \in I$ be an index such that $\frac{1}{A}B_i:=B(x_i,d(x_i))$ do not
verify \eqref{bon2}. Then there is a MS-competitor $L$ for $K$ in
 $$\tilde{B}_i:=B(x_i, \frac{M}{\bar r}d(x_i))$$ such that
$$H^{2}(K\cap \tilde{B}_i)-H^2(L\cap  \tilde{B}_i )\geq
\eta_0 \tilde{r}_i^2 $$ with $\tilde{r}_i:=\frac{M}{\bar r}d(x_i)$
and $M$ is a constant equal to $1$, $10^5$ or $10^{10}$.
\end{prop}
{\bf Proof :} Since $B_i$ do not verify \eqref{bon2}, we know that
$$\beta(x_i,d(x_i))\geq \varepsilon_0$$
and in addition
$$\beta(x_i,\frac{1}{\bar r}d(x_i))\leq \varepsilon_0$$
Multiplying if necessary the radius by $10^5$ or $10^{10}$, and by
use of the re-centering Lemma \ref{recentrage}  (with constant
$V=10^5$), we can suppose that the center of the cone is in a ball
of radius $10^{-5}$ times smaller in $B(x_i,\frac{M}{\bar
r}d(x_i))$ ($M$ is the constant equal to $1$, $10^5$ or
$10^{10}$). Set
$$\tilde {r_i}:=\frac{M}{\bar
r}d(x_i)$$ Then if $\varepsilon_0$ is small enough compared to
$\eta_2$ we have that
$$\beta(x_i , \tilde {r_i})\leq \varepsilon_0 \leq \eta_2$$
with a cone centered in $B(x_i,10^{-5}\tilde{r_i})$. Moreover we
have
$$\beta(x_i,\bar r \tilde{r_i})\geq \frac{1}{M}\varepsilon_0.$$
We also have $F\cap B(x_i,\tilde{r_i})$, that is a separating set
in $B(x_i,\tilde{r_i})$ and such that
$$H^2(F\cap B(x_i,\tilde{r_i}))-H^2(K\cap B(x_i,\tilde{r_i}))\leq  \varepsilon_0'
\tilde{r_i}.$$

 Therefore, we can apply lemma \ref{compacite} in $B(x_i,\tilde{r_i})$
with $\frac{1}{M}\varepsilon_0$ instead of $\varepsilon_0$ that we
may suppose smaller than $C\varepsilon_1$. We can also take
$\varepsilon_0'<<\eta_0$.
 Finally, Lemma \ref{compacite} is stated in $B(0,1)$ but by translation and dilatation it stays true in
 every ball  $B(x,r)$. \qed


\begin{rem}({\rm Choice of $A$}) We can now fix our constant $A$. We
want that for every bad ball $B_i:=B(x_i,Ad(x_i))$ with $i \in I$,
the ball
 $$B(x_i,\tilde{r_i}):=B(x_i,\frac{M}{\bar r}d(x_i))\subset B(x_i,\frac{10A}{U}d(x_i)) \subset \mathpzc{V}$$
in order to have that the extension of $u$ given by  Lemma
\ref{extentionW} is well defined in each $B(x_i,\tilde{r_i})$.
Thus it suffices to take for instance
$$A=\frac{U10^{-20}}{\bar r}.$$
\end{rem}

Before continuing, it is time now to recapitulate in which order
the principal constants are introduced, to see who is controlled
by who. Recall that at beginning we have a Mumford-Shah minimizer
$K$ with  $\beta(0,r_0)$ less than a certain $\varepsilon$. Then
we use a stopping time argument about being close to cones at
small scales with stopping constant $\varepsilon_0$ for the
geometry and $\varepsilon_0'$ for the topology (separating
condition). We obtain a collection of balls that we call ``small
scales'' on which we do some manipulations.

At small scales : The regularity theorem of Guy David gives a
$\varepsilon_1$ for which $\beta$ decays like a power of radius
for a minimal set with excess density (function $f(0,r_0)$)
smaller than $\varepsilon_1$. An other lemma controls $f(0,r_0)$
by $\beta(0,r_0)$ whenever $\beta(0,r_0)$ is smaller than a
certain $\eta_1$. Thus we obtain $\bar r$, that depends on
$\varepsilon_0$, for which $\beta(0,\bar
r)<\frac{\varepsilon_0}{2}$ for all minimal set that is separating
in $B(0,1)$ and such that $\beta(0,1)<\eta_2$. In the proof of
this compactness lemma we fix $\varepsilon_0$ small enough
compared to $\bar \eta_2$. The lemma gives a $\eta_0$ that is the
winning of surface in each bad ball, depending on $\varepsilon_0$
and $\bar r$. In the other hand we have to be sure that
$\varepsilon_0'$ is smaller than
 $\eta_0$ to apply the Lemma in future.
So at this stage we have (each quantity is depending on what is on
the right of the symbol  $\prec$) :
\begin{eqnarray}
\varepsilon_0' \prec \eta_0 \prec \bar r \prec \varepsilon_0 \prec
\bar \eta_2 \label{tab1}
\end{eqnarray}
At big scale : In the big scale we want to show that some
quantities in the ball of radius $ar_0$ are controlled by the same
quantity in the ball of radius $r_0$, for a certain  $a$ that is
chosen later with some arithmetical conditions, in particular
$a^{\gamma}<\frac{1}{8}$ where $\gamma$ is close as we want to
$0,8$. We apply Theorem 9 of \cite{l3} with $F$ a
$(\varepsilon_0,\frac{\varepsilon}{\varepsilon_0})$-minimal set
and  $\varepsilon_0$ is like in the above paragraph. Theorem 9 of
\cite{l3} gives a $\varepsilon_2$ (depending on $\varepsilon_0$,
$\gamma$ and $a$) and assure a decay of energy if $\varepsilon$ is
small enough in respect with $\varepsilon_2$ and $\varepsilon_0$.
Thus in addition of \eqref{tab1} we have
$$\varepsilon\prec \varepsilon_2 \prec
\left\{
\begin{array}{r}
\varepsilon_0 \prec \varepsilon_1 \\
a,\gamma
\end{array}
\right.
$$

\subsection{Bounds for the bad mass}

 The following proposition is an estimate about $m$. Recall that $\rho$
is the radius chosen in $[\frac{r_0}{2},\frac{3}{4}r_0]$.

\begin{prop}\label{comptage} If $m(0,\frac{\rho}{2})>\frac{m(0,r_0)}{10}$ then
\begin{eqnarray}
 m(0,\frac{\rho}{2}) \leq \frac{C}{\eta_0}\left(\omega_2(0,r_0) + \omega_2(0,r_0)^{\frac{1}{2}}J(0,r_0)^{-1}+
 h(r_0)\right).
 \label{compte}
\end{eqnarray}
\end{prop}

{\bf Proof :}  To prove Proposition \ref{comptage}, we will count
the number of $B_i$ for $i \in I$ and use Proposition \ref{compa}
to say that there are not so many. Recall that the $B_i$ are
disjoints.

In order to estimate the bad mass we will take a good competitor
for $(u,K)$ in $B(0,r_0)$. Set $I_1$ the indices of bad balls
$B_i$ such that $B(x_i,d(x_i))$ don't verify \eqref{bon2} and
$I_2:= I \backslash I_1$. In particular, balls of $I_2$ don't
verify \eqref{bon1}. Hence we know that if $i \in I_2$ and if
$r_i:= d(x_i)$ we have
$$r_i^2\leq \frac{1}{\varepsilon_0'}\left(H^2(F\cap B(x_i,r_i))-H^2(K\cap B(x_i,r_i))\right)$$
and since the $B_i$ are disjoint we deduce that
$$\sum_{i\in I_2}r_i^2\leq C\frac{1}{\varepsilon_0'}(H^2(F(0,r_0))-
H^2(K\cap B(0,r_0)))\leq
Cr_0^2\omega_2(0,r_0)^{\frac{1}{2}}J(0,r_0)^{-1}.$$ Now we have to
count the contribution of $I_1$. We will modify each $B_i$ for $i
\in I_1$ with the use of Proposition \ref{compa}. Set
$$
\tilde{G}:=\left\{
\begin{array}{cc}
F(0,r_0) &\text{ in } B(0,r_0) \backslash \bigcup_{i \in I_1 ; B_i \cap B(0,\rho) \not = \emptyset } B_i \\
L_i & \text{ in } B_{i} \text{ for all } i \in I_1 ; B_i \cap B(0,\rho) \not = \emptyset \\
\end{array}
\right.
$$

where $L_i$ is the set given by Proposition \ref{compa}. Then set
$$G:=\tilde{G}\cup \bigcup_{i \in I_\rho}\partial B_i.$$
For the function we use the extension of Proposition
\ref{extentionW} which can be applied in $B(0,\rho)$ by
Proposition \ref{hypothesisH}. Thus we take
$$
v= v^k  \text{ in }   \Omega^k.
$$
By choice of constant $A$ we know that the function $v$ is well
defined in $B(0,r_0) \backslash G$. Set
$$I_1':=\{i \in I_1 ; B_i \cap
B(0,\rho)\not = \emptyset \text{ and } B_i\cap \partial B(0,\rho)=
\emptyset\}$$ and
$$
I_1'':=\{i\in I_1; B_i \cap
\partial B(0,\rho)\not = \emptyset\}.
$$
Notice that  $m(0,\frac{\rho}{2})\leq C\frac{1}{r_0}\sum_{i \in
I_1'}r_i^2$ and $\sum_{i \in I_1''}r_i^2\leq
\sqrt{\varepsilon}m(0,r_0)$. In addition $G$ is a competitor. To
see this we can use the same argument as Remark 1.8. in \cite{d3}.
We apply now the fact that $(u,K)$ is a Mumford-Shah minimizer and
we obtain
\begin{eqnarray}
\int_{B(0,r_0)\backslash K} |\nabla u|^2 + H^2(K\cap B(0,r_0))& \leq &
\int_{B(0,r_0)\backslash G}|\nabla v|^2+ H^2(G\cap B(0,r_0))+r_0^2h(r_0) \notag
\end{eqnarray}

$$\leq
C\int_{B(0,r_0)\backslash K}|\nabla
u|^2+H^2(F(0,r_0))\displaystyle -\eta_0\displaystyle \sum_{i \in
I_1'} r_i^2+ C \sum_{i \in I_1'' } r_i^2+C\sum_{i \in
I_2}r_j^2+r_0^2 h(r_0).$$

Hence,
$$\eta_{0} Cr_0m(0,\frac{\rho}{2}) - C\sqrt{\varepsilon}r_0m(0,r_0)\leq C\int_{B(0,r_0)\backslash K}|\nabla u|^2
+r_0^2\omega_2(0,r_0)^{\frac{1}{2}}J(0,r_0)^{-1}+  r_{0}^2
h(2r_0).$$ Therefore, if $\varepsilon$ is small enough compared to
$\eta_0$ and since $m(0,\frac{\rho}{2})\geq \frac{0,m(r_0)}{10}$
we deduce
\begin{eqnarray}
 m(0,\frac{\rho}{2}) \leq \frac{C}{\eta_0}\left(\omega_2(0,r_0) + \omega_2(0,r_0)^{\frac{1}{2}}J(0,r_0)^{-1}+ h(2r_0)\right)
\notag
\end{eqnarray}
and the proposition follows. \qed


Now by the same sort of argument as Proposition before, we have
this second estimate about $m$.

\begin{prop}\label{estimm}
\begin{eqnarray}
 m(0,r_0(1-5\sqrt{\varepsilon})) \leq \frac{C}{\eta_0}\left(\omega_2(0,r_0) +
  \omega_2(0,r_0)^{\frac{1}{2}}J(0,r_0)^{-1}+
 \beta(0,r_0)+ h(r_0)\right).
\end{eqnarray}
\end{prop}

{\bf Proof :}  The proof is very similar to Proposition
\ref{comptage}. We modify each $B_i$ for $i \in I_1$ with the use
of Proposition \ref{compa}. Set

$$
\tilde{G}:=\left\{
\begin{array}{cc}
F(0,r_0) &\text{ in } B(0,r_0) \backslash \bigcup_{i \in I_1 } B_i \\
L_i & \text{ in } B_{i} \text{ for all } i \in I_1  \\
\end{array}
\right.
$$

where  the $L_i$ are the sets given by Proposition \ref{compa}.
Our competitor is now
$$G:=\tilde{G}\cup T_{\beta}$$

where $T_\beta$ is a little wall of size $\beta:=10\beta(0,r_0)$
$$T_{\beta}:= \{y \in \partial B(0,r_0); d(y,Z)\leq \beta r_0\}$$
with $Z$ a minimal cone centered at the origin at distance less
than $\beta(0,r_0)$ of $K$ in $B(0,r_0)$.

We keek the same notation $I_1$, $I_2$, $I_1'$ and $I_1''$ as
before but now with $\rho=r_0$. As before we have
$$\sum_{i\in I_2}r_i^2\leq C\frac{1}{\varepsilon_0'}(H^2(F(0,r_0))-H^2(K\cap B(0,r_0)))\leq Cr_0^2\omega_2(0,r_0)^{\frac{1}{2}}J(0,r_0)^{-1}.$$

For the function we use the extension of Proposition
\ref{extentionW} in $B(0,2r_0)$ with $\rho=r_0$ and with
application $\delta$ defined in \eqref{defd}. We set
$$
v= v^k  \text{ in }   \Omega^k.
$$
By choice of constant $A$ we know that the function $v$ is well
defined in $B(0,r_0) \backslash G$ and since we added $T_{\beta}$
there is no boundary problem.

We apply now the fact that $(u,K)$ is a Mumford-Shah minimizer and
we obtain with same notations as Proposition before,
\begin{eqnarray}
\int_{B(0,r_0)\backslash K} |\nabla u|^2 + H^2(K\cap B(0,r_0))&
\leq & \int_{B(0,r_0)\backslash G}|\nabla v|^2+ H^2(G\cap
B(0,r_0))+r_0^2h(r_0) \notag
\end{eqnarray}

$$\leq
C\int_{B(0,r_0)\backslash K}|\nabla u|^2+H^2(F(0,r_0))
 -\eta_0\displaystyle \sum_{i \in I_1' } r_i^2
+C\sum_{i \in I_2}r_j^2+H^2(T_\beta)+r_0^2 h(r_0)$$

Hence,
$$\eta_{0} m(0,r_0(1-5\sqrt{\varepsilon})) \leq C\int_{B(0,r_0)\backslash K}|\nabla u|^2
+r_0^2\omega_2(0,r_0)^{\frac{1}{2}}J(0,r_0)^{-1}+
Cr_0^2\beta(0,r_0)+r_{0}^2 h(2r_0)$$ because all the $B_i$ have a
radius less than $\sqrt{\varepsilon}r_0$ thus all the $B_i$ for $i
\in I_1 $ such that $5B_i \cap
\partial B(0,r_0) = \emptyset$ are included in
$B(0,r_0(1-5\sqrt{\varepsilon}))$,
and the proposition follows. \qed


\subsection{Control of the minimality defect }

In this section we want to control the defect of minimality of $K$
in terms of energy and bad mass. For some topological reasons we
are not going to work directly on $K$, but we will use the set $F$
to be sure that it is separating in $B$. We show in this section
that for all MS-competitor $L$ for $F$ we can give a function  $w$
such that $(L,w)$ is a Mumford-Shah competitor for $(u,K)$, and
with good bounds on the energy of $w$. Here is a more precise
statement:

\begin{prop} \label{comm}  There is a positive constant $c_{10}<1$ such that
for all MS-competitor $L$ for the set $F$ (see Definition
\ref{m1})
in the ball $B(0, c_{10}r_0)$, we have :\\
$$\displaystyle{\frac{1}{r_0^2}[H^{2}(F\cap
B(0,c_{10}r_0))-H^{2}(L\cap B(0,c_{10}r_0))]} \leq
\displaystyle{C\left[\omega_2(0,r_0)+\sqrt{\varepsilon}m(0,r_0)+h(r_0)\right]}
$$
\end{prop}

{\bf Proof :} Let $Z^0$ be the cone such that $d(x,Z^0)\leq
\varepsilon r_0$ for all $x \in K \cap B(0,r_0)$. We call as
usuall $Z^0_\varepsilon$ the region
\begin{eqnarray}
Z^0_\varepsilon := \{x \in B(0,r_0); d(x,Z_0) \leq \varepsilon \}.
\label{defTlambda}
\end{eqnarray}
We consider our ball $\{B_i\}_{i \in I}$ obtained by the stopping
time argument. We define the functions
$$\psi_i:=\left\{
\begin{array}{cc}
r_i & \text{ on } B_i\\
0 &\text{ in the complement of } 2B_i
\end{array}
\right.
$$

then for all $x$ we define
$$d_1(x):=\sum_{i \in I} \psi_{i}(x).$$
Finally, for all $x \in B(0,\rho)$ set
$$\delta(x):= \max(d(x, \partial B(0,\rho)),d_1(x)).$$

As usual, $\delta(x)$ is a geometric function associated to $F$ in
$B(0,r_0)$. Thus applying Lemma \ref{extentionW} we get
$\mathpzc{k}^\rho$ functions $v^k$ such that  $ v^{k} \in
W^{1,2}(\Omega^k \cup \mathpzc{V})$ and such that
$$\int_{\Omega^k \cup \mathpzc{V}\backslash \mathpzc{V}_\rho}|\nabla v^{k}|^2 \leq C\int_{B(0,\rho)\backslash F}|\nabla u|^2$$
in addition,  $ v^{k}$ is equal to  $u$ on $\partial B(0,\rho)
\cap \Omega^{k}\backslash \mathpzc{V}$.

Moreover, since $\delta(x)\geq d(x,\partial B(x, \rho))$, if
$\varepsilon$ is small enough we can easily deduce that there is a
constant  $c_{10}<\frac{1}{2}$ depending on constant $U$ such that
$B(0,c_{10}r_0)\subset \mathpzc{V}$. Set
$$G'= \left\{
\begin{array}{cc}
F & \text{ in } B(0,r_0)\backslash B(0,c_{10}r_0) \\
L & \text{ in } B(0,c_{10}r_0)
\end{array}\right.
$$
If $L$ is a MS-competitor for $F$ in $B(0,c_{10} r_0)$, we know
that $L$ in separating $B(0,c_{10}r_0)$ into $\mathpzc{k}^\rho$
big connected components (because $F$ is separating and $L$ is a
topological competitor). Thus $G'$ is separating in $B(0,\rho)$
ant we note $(B(0,\rho)\backslash G)^k$ the big connected
components.

Then set

$$G:=G' \cup \bigcup_{i \in I(\rho)}\partial B_i$$

and

$$
v:=\left\{
\begin{array}{ccc}
u & \text{ in } & B(0,r_0)\backslash B(0,\rho) \\
v^{k} & \text{ in } &  (B(0,\rho)\backslash G)^k\\
0 &\text{ in other components of } &  B(0,\rho)\backslash G
\end{array}
\right.
$$

Using that $(u,K)$ is a  Mumford-Shah minimizer and that $(v,G)$
is a competitor we obtain
\begin{eqnarray}
\int_{B(0, \rho) \backslash K}|\nabla u|^2 + H^2(K) &\leq&
\int_{B(0,\rho) \backslash G}|\nabla v|^2 + H^2(G)+\rho^2 h(\rho)
\notag
\end{eqnarray}

thus

$ \displaystyle{H^{2}(K\cap B(0,c_{10} r_0))-H^2(L\cap
B(0,c_{10}r_0))} $
\begin{flushright}
$ \displaystyle{\leq C\left[\int_{B(0, r_0) \backslash K}|\nabla
u|^2   + \sum_{i\in I_1''}r_i^2
+r_0^{2}\omega(0,r)^{\frac{1}{2}}J(0,r)^{-1} +r_0^2 h(r_0)\right]}
$
\end{flushright}
and the proposition follows. \qed

\subsection{Conclusion about regularity}

Now we are ready to use all the preceding estimates in order to
prove some regularity. We begin with this proposition about
self-improving estimates.

\begin{prop} \label{self} There is an $\varepsilon>0$, some   $\tau_4<\tau_3<\tau_2<\tau_1<\varepsilon$ and  $a<1$
such that if $x \in K$ and $r$ are such that $B(x,r)\subset
\Omega$, and
\begin{eqnarray}
h(r)+J(x,r)^{-1}\leq \tau_4, \quad  \omega_2(x,r)\leq \tau_3,
\quad m(x,r)\leq \tau_2 , \quad \beta(x,r)\leq \tau_1 \label{tau4}
\end{eqnarray}
then \eqref{tau4} is still true with $ar$ instead of  $r$.
\end{prop}
{\bf Proof :} We choose  $\varepsilon<<\varepsilon_0$ and
$\varepsilon_1$ such that all the results of the preceding
sections are true. We choose $a<\frac{1}{16}$ such that applying
\eqref{debdec} to $(u,K)$ gives
\begin{eqnarray}
\omega_2(x,ar)\leq
\frac{1}{8}\omega_2(x,r)+C_2\omega_2(x,r)^{\frac{1}{2}}J(x,r)^{-1}+C_2\sqrt{\varepsilon}m(x,r)+C_2h(r).
\label{goq}
\end{eqnarray}
Since $a$ is chosen, we can fix  $\tau_1$ small enough such that
for all $ar<t<r$ we have $\beta(x,t)\leq 10^{-1}$. Hence by Lemma
\ref{lemme7}
$$J(x,ar)  \geq a^{-\frac{1}{2}}[J(x,r)-C']\geq \frac{1}{2}a^{-\frac{1}{2}}J(x,r)$$
if $\tau_4$ is small enough compared to $C'$. Then we deduce
$$J(x,ar)^{-1}\leq 2a^{\frac{1}{2}} J(x,r)^{-1}\leq \frac{J(x,r)^{-1}}{2}$$
because $a<\frac{1}{16}$. In addition if $\tau_4$ is small enough
compared to $\tau_3$, we have
\begin{eqnarray}
C\tau_3^{\frac{1}{2}}\tau_4\leq\frac{1}{8} \tau_3 . \label{cond0}
\end{eqnarray}
Therefore by \eqref{goq},
$$\omega_2(x,ar)\leq \frac{3}{8}\tau_3+C_2\sqrt{\varepsilon}m(x,r)\leq \frac{\tau_3}{2}$$
under the condition that
\begin{eqnarray}
8C_2\sqrt{\varepsilon}\tau_2 <\tau_3 .\label{cond1}
\end{eqnarray}

Now for $m(x,r)$ we have two cases. If $m(x,ar)\leq
\frac{m(x,r)}{10}$ then $m(x,ar)\leq \frac{\tau_2}{10}$ and it is
what we want. Otherwise, we have $m(x,ar)>\frac{m(x,r)}{10}$ which
implies $m(x,\frac{\rho}{2})>\frac{a^2 m(x,r)}{5}$ and then we can
use the proof of Proposition  \ref{comptage} with a slightly
different constant (depending on $a$) to obtain
$$m(x,ar)\leq \frac{C(a)}{\varepsilon_1}(\tau_3+\tau_3^{\frac{1}{2}}\tau_4+\tau_4)\leq
\frac{C(a)}{\varepsilon_1}\tau_3 \leq \frac{\tau_2}{2}$$ if
\begin{eqnarray}
2\frac{C(a)}{\varepsilon_1}\tau_3\leq \tau_2 .\label{cond2}
\end{eqnarray}
So it suffice to choose  $\varepsilon$ small enough compared to
 $\varepsilon_0$ and $C$ in order to have the existence of
$\tau_3<\tau_2$ that verify simultaneously \eqref{cond1} and
\eqref{cond2}. Hence, we control $\omega_2(x,ar)$ and $m(x,ar)$. \\

To finish we have to control  $\beta(x,ar)$. For that we use the
estimate in Proposition  \ref{comm} and  Lemma \ref{compacite}
that we apply in $B(x,c_{10}r)$. Indeed, suppose that $a<<\bar
r(\varepsilon_0)$ is such that
\begin{eqnarray}
\beta(x,ar)\geq \tau_1 \label{beta5}
\end{eqnarray}
 Then applying Lemma \ref{compacite} with
$\varepsilon_0=\tau_1$ gives a $\eta_0(\tau_1,a)$ and a competitor
$L$ for $K$ in $B(x,c_{10}r)$ such that
\begin{eqnarray}
H^2(K)-H^2(L)\geq \eta_0(\tau_1,a) \label{in}.
\end{eqnarray}
On the other hand, according to Proposition $\ref{comm}$, if we
choose $\tau_2$ and $\tau_3$ small enough compared to
$\eta_0(\tau_1,a)$, the inequality \eqref{in} cannot hold. This
shows that
$$\beta(x,ar)\leq \tau_1$$
and  gives a contradiction with \eqref{beta5} which achieves the
proof of the proposition. \qed

We keep the constants  $a$ and $\tau_i$ given by the preceding
proposition. Let $b$ be the positive power such that
$a^{b}=\frac{1}{2}$. Set
$$\tilde h_r(t)=\sup\big\{\left(\frac{t}{s}\right)^{b}h(s) ; t\leq s \leq r\big \}$$
for $t<r$ and $\tilde h_r(t)=h(t)$ for $t>r$. According to
\cite{d} page 318, the function $\tilde h$ is still a gauge
function (i.e. monotone and with limit equal to $0$ at $0$). We
also trivially have that $h(t)\leq \tilde h_r(t)$ and one can
prove that
\begin{eqnarray}
\tilde h_r(t)\geq \left(\frac{t}{t'}\right)^b\tilde h_r(t') \quad
\text{ for }0<t<t'\leq r.
\end{eqnarray}

Note that since $a^b=\frac{1}{2}$, we have
\begin{eqnarray}
\tilde h_r(at)\geq \frac{1}{2}\tilde h_r(t) \quad \text{ for } 0<
t \leq r. \label{truc}
\end{eqnarray}

The purpose of Proposition \ref{self} is just to have
$\beta(x,r)\leq \tau_1$ at all scales in order to have more decay
for the other quantities. Notice that at this step, we could prove
that $K$ is the bi-hölderian image of a minimal cone using
\cite{ddpt}. This will be done in Corollary \ref{fegk} to prove
that $K$ is a separating set. Before that we will prove some more
decay estimates.

\begin{prop} We assume that we have the same hypothesis as in the proposition before.
Then for all $0<t<r$ we have
$$J(x,t)^{-1}\leq 2\left(\frac{t}{r}\right)^{b}\tau_4$$
$$\omega_2(x,t)\leq C\left(\frac{t}{r}\right)^{b}\tau_3 + C\tilde h_r(t)$$
$$m(x,t) \leq C\left(\frac{t}{r}\right)^b\tau_2+C\tilde h_r(t).$$
\end{prop}
{ \bf Proof :} The first step is to control the jump. Since
$\tau_1$ is small enough to have $\beta(x,t) \leq 10^{-1}$ for all
$t<r$, then by Lemma \ref{lemme7}
$$J(x,t)  \geq \left(\frac{r}{t}\right)^{-\frac{1}{2}}[J(x,r)-C']\geq \frac{1}{2}
\left(\frac{r}{t}\right)^{-\frac{1}{2}}J(x,r)$$ if $\tau_4$ is
small enough compared to $C'$. We deduce
\begin{eqnarray}
J(x,t)^{-1}\leq 2\left(\frac{t}{r}\right)^{\frac{1}{2}}
J(x,r)^{-1}. \notag
\end{eqnarray}
And since $a<\frac{1}{4}$ we have
\begin{eqnarray}
J(x,a^nr)^{-1}\leq 2\left(\frac{1}{2}\right)^{n}\tau_4.
\label{decrsaut}
\end{eqnarray}
 Now we want to show by induction
that
\begin{eqnarray}
\omega_2(x,a^nr_0)\leq 2^{-n}\tau_3+C_3\tilde h_r(a^nr) \quad
\text{ and } \quad  m(x,a^nr)\leq 2^{-n}\tau_2+ C_3\tilde
h_r(a^nr) \label{rec}
\end{eqnarray}

For $n=0$ we have \eqref{rec} trivially. Suppose by now that
\eqref{rec} is true for $n$.
 Then applying inequality  \eqref{debdec}  in $B(x,a^nr)$
\begin{eqnarray}
\omega_2(x,a^{n+1}r)\leq \frac{1}{8}\omega_2(x,a^nr)+C_2\omega_2
(x,a^nr)^{\frac{1}{2}}J(x,a^nr)^{-1}+C_2\sqrt{\varepsilon}m(x,a^nr)+C_2h(a^nr).
\label{toto}
\end{eqnarray}
Now, using the inequality $2ab\leq a^2+b^2$ we obtain
$$\omega_2(x,a^{n}r)^{\frac{1}{2}}\leq \frac{1}{20C_2}\omega_2(x,a^nr)J(x,a^nr)+5C_2J(x,a^nr)^{-1}.$$
Thus \eqref{toto} yields
$$\omega_2(x,a^{n+1}r)\leq \frac{7}{40}\omega_2(x,a^nr)+5C_2^2J(x,a^nr)^{-2}+
C_2\sqrt{\varepsilon}m(x,a^nr)+C_2h(a^nr).$$ Now using
\eqref{decrsaut}, and the induction hypothesis we obtain
$$\omega_2(x,a^{n+1}r)\leq \frac{7}{40}2^{-n}\tau_3+5C_2^24\tau_4^22^{-n}+
C_2\sqrt{\varepsilon}2^{-n}\tau_2+(\frac{7}{40}C_3+C_2\sqrt{\varepsilon}C_3+C_2)\tilde
h_r(a^nr).$$ Now, using that $\tau_4$ controlled by $\tau_3$,
since $\varepsilon$ is small as we want compared to $C_2$, using
also \eqref{cond1} and \eqref{truc}, and finally if we choose
$C_3$ larger than $100C_2$ we deduce that
$$\omega_2(x,a^{n+1}r)\leq (\frac{8}{40} +\frac{1}{8})2^{-n}\tau_3+
C_3\tilde h_r(a^nr)\leq 2^{-(n+1)}\tau_3+ C_3\tilde
h_r(a^{n+1}r).$$

Concerning $m(x,r)$ it is a similar argument, suppose that
$m(x,a^{n+1}r)> 2^{-(n+1)} m(x,a^{n+1}r)$. Then we can apply
Proposition \ref{comptage} in the ball $B(x,a^{n}r)$ thus
\begin{eqnarray}
m(x,a^{n+1}r)&\leq&
\frac{C(a)}{\varepsilon_1}(\omega_2(x,a^nr)+\omega_2(x,a^nr)^{\frac{1}{2}}J(x,a^nr)^{-1}+h(a^nr))\notag
\\
 &\leq
 &\frac{C(a)}{\varepsilon_1}(\frac{3}{2}\omega_2(x,a^nr)+\frac{1}{2}J(x,a^nr)^{-2}+h(a^nr)).
\label{deb}
\end{eqnarray}
Setting $C_4=\frac{C(a)}{\varepsilon_1}$, using \eqref{decrsaut}
and induction hypothesis we obtain
\begin{eqnarray}
m(x,a^{n+1}r)&\leq& C_4 2^{-n}\tau_3+C_42^{-n}\tau_4+2C_4\tilde
h_r(a^nr) \notag \\
&\leq & 2^{-n}\tau_2+C_3\tilde h_r(a^{n+1}r) \notag
\end{eqnarray}
because $\tau_3$ and $\tau_4$ are small as we want with respect to
$C_4$ and $\tau_2$, and because we can chose $C_3$ bigger than
$10C_4$ and we have used \eqref{truc}.

To finish the proof let $0<t<r$ and $n$ such that $a^{n+1}\leq t
\leq a^nr$. Then we have
\begin{eqnarray}
\omega_2(x,t)=\frac{1}{t^2}\int_{B(0,t)\backslash K}|\nabla u|^2 &\leq& \left(\frac{a^nr}{t}\right)^2\omega_2(x,a^nr) \notag\\
&\leq& \frac{1}{a^2}2^{-n}\tau_3 +C_3\tilde h_r(a^nr) \notag \\
&\leq& \frac{1}{a^2}a^{bn}\tau_3 +C_3'\tilde h_r(t)  \notag \\
&\leq& C\left(\frac{t}{r}\right)^b \tau_3+C_3'\tilde h_r(t)
\notag
\end{eqnarray}
and
\begin{eqnarray}
m(x,t)&\leq& \frac{a^{2n}r^2}{t^2}m(x,a^{n}r) \leq \frac{a^{2n}r}{t^2} 2^{-n}\tau_2 + C_3\tilde h_r(a^nr) \notag \\
&\leq & \frac{1}{a^2}a^{bn}\tau_2 + C_3'\tilde h_r(t) \leq
C\left(\frac{t}{r} \right)^{b}\tau_2 +   C_3'\tilde h_r(t) .\notag
\qed
\end{eqnarray}

\begin{prop} \label{decryep} There is a positive constant $b$ such that the following is true. Let $(u,K)$
be a Mumford-Shah minimizer in $\Omega \subset \R^3$ with gauge
function $h$. Let  $x_0 \in K$ and $r_0$ be such that
$B(x_0,r_0)\subset \Omega$.
 Then there is
$\varepsilon>0$ and $\tau_4'<\tau_3'<\tau_2'<\tau_1'<\varepsilon$
such that if
\begin{eqnarray}
h(r_0)+J(x_0,r_0)^{-1}\leq \tau_4', \quad  \omega_2(x_0,r_0)\leq \tau_3',
\quad  m(x_0,r_0)\leq \tau_2', \quad \beta(x_0,r_0)\leq \tau_1' \notag
\end{eqnarray}
then for all $x \in B(x_0,\frac{1}{10}r_0)$ and for all
$0<t<\frac{1}{2}r_0$ we have
$$J(x,t)^{-1}\leq C\left(\frac{t}{r_0}\right)^{b}$$
$$\omega_2(x,r)\leq C\left(\frac{t}{r_0}\right)^{b} + C\tilde h_r(t)$$
$$m(x,t) \leq C\left(\frac{t}{r_0}\right)^b+C\tilde h_r(t)$$
$$\beta(x,t) \leq \tau_1$$
\end{prop}
{\bf Proof :} It suffice to show that there is
$\tau_4'<\tau_3'<\tau_2'<\tau_1'<\varepsilon$ such that if
\begin{eqnarray}
h(r_0)+J(x_0,r_0)^{-1}\leq \tau_4', \quad  \omega_2(x_0,r_0)\leq \tau_3',
\quad  m(x_0,r_0)\leq \tau_2', \quad \beta(x_0,r_0)\leq \tau_1' \notag
\end{eqnarray}
then for all $x \in B(x_0,\frac{1}{10}r_0)$ we have
\begin{eqnarray}
h(\frac{1}{2}r_0)+J(x,\frac{1}{2}r_0)^{-1}\leq \tau_4, \quad
\omega_2(x,\frac{1}{2}r_0)\leq \tau_3, \quad  m(x,\frac{1}{2}r_0)\leq
\tau_2, \quad \beta(x,\frac{1}{2}r_0)\leq \tau_1 \notag
\end{eqnarray}
hence we could apply all the work of preceding sections in
$B(x,\frac{1}{10}r_0)$ and conclude.

Note that for all $x \in K\cap B(x_0,\frac{1}{10}r_0)$ we have
\begin{eqnarray}
\omega_2(x,\frac{1}{2}r_0) \leq 4\omega_2(x_0,r_0) \notag \\
m(x,\frac{1}{2}r_0) \leq 2 m(x_0,r_0) \notag \\
\beta(x,\frac{1}{2}r_0) \leq 2\beta(x_0, r_0)
\end{eqnarray}
in addition if $\beta(x_0,r_0) $ is small enough then
$$J^{-1}(x,\frac{1}{2}r_0)\leq 2 J^{-1}(x_0,r_0).$$
Finally, since $h$ is non decreasing
$$h(\frac{1}{2}r_0)\leq h(r_0).$$
We deduce that for  $i \in [1,4]$ we can set
$$\tau_i':=\frac{1}{4}\tau_i$$
and the proposition follows. \qed

\begin{cor}\label{fegk} In the same situation as in proposition before,
if $\tau_1$ is small enough we can choose
$$F(x_0,\frac{1}{10}r_0)=K\cap B(x_0,\frac{1}{10}r_0).$$
\end{cor}
{\bf Proof :} The method is to prove that $K$ is separating in
 $B(x_0,\frac{1}{10}r_0)$. This will show that we can take $F=K$ in this ball.
 To show that $K$ is separating we will apply Theorem 1.1 of \cite{ddpt}, even if we could prove
 the same result without
 using  \cite{ddpt} but with a longer explication. The main point is to show that
 for all $x \in B(x_0,\frac{1}{10}r_0)$ and for all $r$ such that
$B(x,r)\subset B(x_0,\frac{1}{5}r_0)$ there is a cone $Z(x,r)$
such that
$$D_{x,r}(K,P(x,r))\leq \varepsilon'$$
with $\varepsilon'$ a certain constant given by Theorem 1.1. of
\cite{ddpt}. Recall that according to the notations of
\cite{ddpt}, $D_{x,r}$ is the Hausdorff distance
\begin{eqnarray}
D_{x,r}(E,F):=\frac{1}{r}max\big\{\sup_{z \in E\cap
B(x,r)}\{d(z,F)\},\sup_{z \in F\cap B(x,r)}\{d(z,E)\}\big\}.
\label{dxr}
\end{eqnarray}
 If we
choose  $\tau_1$ small enough compared to $\varepsilon'$ we know
that for all $x$ and for all $r$ we have $\beta(x,r)\leq
\varepsilon'$ by the preceding proposition. Hence we can find a
cone $Z(x,r)$ that satisfy the first half of  $D(x,r)$. We have to
show now that
$$\sup \{d(z,K), z \in Z(x,r)\}\leq r\varepsilon'.$$
We know that $J(x,r)^{-1} \leq \tau_4$ and $\omega_2(x,r)\leq
\tau_3$. Thus there is a set $F(x,r)$ that is separating in
$B(x,r)$ and such that
$$H^2(F(x,r)\cap K\cap B(0,r))\leq C\omega_2(x,r)^{\frac{1}{2}}J(x,r)^{-1}\leq \frac{1}{8}\tau_3r^2.$$
Then for all $z \in Z(x,r)$, we have
$$d(z,K)\leq d(z,y)+d(y,K)$$
with $y$ a point of $F(x,r)$ such that $d(z,F(x,r))=d(z,y)$. If
$\tau_1 \leq \frac{\varepsilon'}{2}$ we can suppose that $F(x,r)
\subset \{y ; d(y,Z) \leq r\frac{\varepsilon'}{2}\}$. Thus $d(z,y)
\leq r\frac{\varepsilon'}{2}$. We claim that  $d(y,K) \leq
r\frac{\varepsilon'}{2}$. The argument is by contradiction. If it
is not true, then $K\cap B(y,r\frac{\varepsilon'}{2})=\emptyset$.
But $F(x,r)$ is included in $T:=\{y ; d(y,Z) \leq
r\varepsilon'\}$. Let $A^k$ be the connected components of
$B(y,r\frac{\varepsilon'}{2})\backslash T$. Then $F(x,r)$
separates the $A^k$ in $B(y,r\frac{\varepsilon'}{2})$, and the
minimal set that have this property is a cone of type $\Pp$, $\Y$
or $\T$ of area greater than $C\varepsilon'^2 r^2$. On the other
hand $H^2(F(x,r)\backslash K)\leq \tau_3r^2$. Thus if $\tau_3$ is
small enough compared to $\varepsilon'$ it is not possible, thus
finally $d(y,K)\leq \frac{\varepsilon}{2}$ and
$$D_{x,r}(K,P)\leq \varepsilon'.$$
Now  Theorem 1.1 of \cite{ddpt} says that $K$ is containing the
image of a minimal cone by a homeomorphism from
$B(x_0,\frac{1}{10}r_0)$ to $B(x_0,\frac{1}{5}r_0)$. This proves
that $K$ separates $D^+$ from $D^-$ in
$B(x_0,\frac{1}{10}r_0)$.\qed

\begin{theo} \label{fin2} There is some absolute positive constants  $\varepsilon$ and $c$ such that the following is
true. Let $(u,K)$ be a Mumford-Shah minimizer  in $\Omega \subset
\R^3$ with gauge function $h$, let $x \in K$ and $r$ be such that
$B(x,r)\subset{\Omega}$ and
$$\omega_2(x,r)+ \beta(x,r)+J(x,r)^{-1}+h(r)\leq \varepsilon$$
where the best cone in $\beta(x,r)$ named $Z$ is of type $\Pp$,
$\Y$ or $\T$ centered at $x$. Then there is a diffeomorphism
$\phi$ of class $C^{1,\alpha}$ from $B(x, cr)$ to its image such
that $K \cap B(x, cr)=\phi(Z)\cap B(x,cr)$.
\end{theo}
{\bf Proof :} We want to apply Corollary 12.25 of \cite{d5} (or
see Corollary \ref{gd}).

Thus to prove Theorem  \ref{fin2}, it suffice to show that $K\cap
B(x,cr)$ is an almost minimal set  that verify tht hypothesis of
Corollary \ref{gd}. If  $\varepsilon$ is small enough,
 all the quantities $\omega_2(x,r)$, $\beta(x,r)$, $J^{-1}(x,r)$ and $h(r)$ verify the hypothesis of
 Proposition \ref{decryep}.
   In addition, according to Proposition \ref{estimm}
   (applied in $B(x,r(1-\sqrt{\varepsilon})^{-1})$), $m(r)$ is also smaller that $\tau_2$. So we can apply the
   result of the preceding Propositions.

By Corollary  \ref{fegk}, we know that $F=K$ in $
B(x,\frac{1}{10}r)$. So we can apply Proposition \ref{comm}
directly on $K$ (instead of  $F$) and the monotonicity of $\omega$
and $m$ obtained in Proposition \ref{decryep} shows that $K$ is an
almost minimal set in $B(x,\frac{1}{10}c_{10}r)$ with gauge
function
$$\hat h(t):=C\left(\frac{t}{r}\right)^b+C\tilde h_r(t).$$
To conclude we have to verify \eqref{defautdens}. If $\varepsilon$
and $c$ are small enough we have that $\hat h(cr)\leq
\varepsilon_1$ so we only have to control $f(x,r)$. To do this we
can use the same argument as we used in Lemma \ref{compa2}. We use
Lemma \ref{compa1} to find a point $x$ of same type of cone $Z$
that define $f$, then we use the same competitor $L$ as in the
proof of \ref{compa2} that is $Z\cup M$ where $M$ is a small wall.
We deduce a bound of $f$ by $\beta$. Thus if the $\tau_i$ are
small enough compared to $\varepsilon_1$, \eqref{defautdens} is
verified hence the proof is achieved. \qed

\begin{rem} Constant $c$ in Theorem \ref{fin2} is depending on
$c_{10}$, $U$, $\alpha$, and other constants. Thus, constant $c$
is fairly small but one might give an explicit value by doing some
long computations.
\end{rem}

Now we want to prove that the conditions on $J$ and $\omega_2$ can
be removed in Theorem \ref{fin2} if we suppose that $c$ and
$\varepsilon$ are a bit smaller. To begin, we have to use this
following lemma.

\begin{lem}\label{jump} There is some absolute positive constants
$\varepsilon_3$ and  $\eta_1$ such that if $x \in K$,
$B(x,r)\subset \Omega$,
$$\omega_2(x,r)+h(r)+\beta(x,r)\leq \varepsilon_3$$
then $J(x,r)\geq \eta_1$.
\end{lem}

{\bf Proof :}  The proof is like Lemma 8 page 365 and Proposition
10 page 297 of \cite{d}. The generalization of these lemmas in
higher dimension is not a problem by the same way as we have
proved Lemma \ref{constrF}, Lemma \ref{lemme1} and Lemma
\ref{lemme7}. \qed

About the normalized energy we also have this result that
naturally comes from an argument with blow up limits. One can find
a similar statement about dimension 2 in Lemma 3 page 504 of
\cite{d}. The proof is the same for the case of $\Y$ and $\T$ in
$\R^3$ so it has been omitted here. Recall that $D_{x,r}$ is the
normalized bilateral Hausdorff distance defined in \eqref{dxr}.

\begin{lem}\label{controlen} For each $\eta_2>0$ there is constants
$\varepsilon_3$ and  $a_0$ with the following property. Let
$\Omega \subset \R^3$ and let $(u,K)$ be a Mumford-Shah minimizer
in $\Omega$ with gauge function $h$. Let $x \in K$ and $r>0$ be
such that $B(x,r)\subset \Omega$. Suppose that $h(r)\leq
\varepsilon_3$ and that we can find a cone $Z$ of type $\Pp$, $\Y$
or $\T$ centered at $x$ such that
$$D_{x,r}(K,Z)\leq \varepsilon_3.$$
Then
$$\omega_2(x,a_0r)\leq \eta_2.$$
\end{lem}

Now we can state the main theorem.

\begin{theo} \label{fin3} There is some absolute positive constants
$\varepsilon$ and $c$ such that the following is true. Let $(u,K)$
be a Mumford-Shah minimizer  in $\Omega \subset \R^3$ with gauge
function $h$, let $x \in K$ and $r$ be such that
$B(x,r)\subset{\Omega}$ and $h(r)\leq \varepsilon$. Assume in
addition that there is a cone $Z$ of type $\Pp$, $\Y$ or $\T$
centered at $x$ such that
$$D_{x,r}(K,Z)\leq \varepsilon.$$
Then there is a diffeomorphism $\phi$ of class $C^{1,\alpha}$ from
$B(x, cr)$ to its image, such that  $K \cap B(x, cr)=\phi(Z)$.
\end{theo}

{\bf Proof :} We have to control the normalized jump and then
apply  Theorem \ref{fin2}. Firstly, if $\varepsilon$ is small
enough compared to $\varepsilon_3$ we can use Lemma \ref{jump} and
obtain that
$$J(x,r)\geq \eta_1$$
for a certain $\eta_1>0$. Then, by Lemma \ref{lemme7} we have, for
$r'\leq r$,
$$J(x,r')\geq  \left(\frac{r}{r'}\right)^{\frac{1}{2}}[J(x,r)-\omega_{2}(x,r)].$$
If $\varepsilon$ is small enough compared to $\eta_1$,  the
quantity $J(x,r)-\omega_{2}(x,r)$ is positive. Then by a good
choice of $r'$, and if $\varepsilon$ is small enough compared to
$\frac{r}{r'}$, we deduce that
$$J(x,r')^{-1}\leq \bar \varepsilon$$
where $\bar \varepsilon$ is the constant of Theorem \ref{fin2}.

Now since $\varepsilon$ is still small as we want, we can assume
that the cone in $\beta(x,r')$ is still centered near $x$ and in
addition
$$\beta(x,r')+ J(x,r')^{-1}+\omega_2(x,r')+h(r')\leq \bar \varepsilon.$$
Then we apply Theorem \ref{fin2} in $B(x,r')$ and the conclusion
follows. \qed

This is an example of statement in terms of functional $J$.

\begin{cor} \label{fin4} There is some absolute positive constants
$\varepsilon$ and $c$ such that the following is true. Let $g \in
L^\infty$ and  $\Omega \subset \R^3$. There is a $\tilde r$ that
depends only on $\|g\|_{\infty}$ such that for all pair $(u,K)\in
\mathcal{A}$ that minimize the functional
$$J(u,K):= \int_{\Omega \backslash K}|\nabla u|^2dx+\int_{\Omega\backslash K}(u-g)^2dx+H^{1}(K),$$
for all $x \in K$ and $r<\tilde r$ such that there is a cone $Z$
of type $\Pp$, $\Y$ or $\T$ centered at $x$ with
$$D_{x,r}(K,Z)\leq \varepsilon$$
 there is a diffeomorphism $\phi$ of
class $C^{1,\alpha}$ from $B(x, cr)$ to $B(x,10cr)$ such that  $K
\cap B(x, cr)=\phi(Z)\cap B(x,cr)$.
\end{cor}
{\bf Proof :} We know by Proposition 7.8. p 46 of \cite{d} that
$(u,K)$ is a Mumford-Shah minimizer with gauge function
$$h(r)=C_N \|g\|_{\infty}^2r$$
where $C_N$ depends only on dimension. The conclusion follows
applying Theorem \ref{fin3} in $B(x,r)$ if we choose
$$\tilde r = \frac{\tilde \varepsilon}{2C_N\|g\|_{\infty}^2}$$
where $\tilde \varepsilon$ is the constant of Theorem \ref{fin3}.
\qed

Now we want a statement with only a condition about energy. We
begin by this following lemma ($D_H$ denotes the Hausdorff
distance).

\begin{lem} \label{compa3} For every $\eta_4>0$ there exist a radius  $R>1$ and a  $\eta_3>0$
 such that for every Mumford-Shah minimizer $(u,K)$ in
$B(x,R)\subset \R^3$ such that $x \in K$ and
$$\omega_2(x,R)+h(R)\leq \delta_3,$$ there is a minimal cone
$Z$ of type $\Pp$, $\Y$ or $\T$ that contains $x$ and such that
$$D_H(K\cap B(0,1),Z\cap B(0,1))\leq \delta_4.$$
\end{lem}
{\bf Proof :} The argument is by compactness. If it is not true,
then we can find a $\eta_4>0$ such that for all $n>0$, there is a
Mumford-Shah minimizer $(u_n,K_n)$ in $B(x,n)$ such that
\begin{eqnarray}
\omega_2(x,n)+h(n)\leq \frac{1}{n^3} \label{lim}
\end{eqnarray}
 and
\begin{eqnarray}
\sup_{Z}D_H(K_n\cap B(0,1),Z\cap B(0,1))\geq \eta_4 \label{derv}
\end{eqnarray} where
the supremum is taken over all minimal cones containing $x$. We
let now tend $n$ to infinity. Since $(u_n, K_n)$ is a sequence of
Mumford-Shah minimizers, with same gauge function
$h_l(r):=\sup\{h(nr); n\geq l\}$, and such that
$$\int_{B(x,n)}|\nabla u|^2 \leq r\frac{1}{n}\leq C$$
by Proposition 37.8 of \cite{d} we can extract a subsequence such
that  $(u_{n_k},K_{n_k})$ converges to $(u,K)$ in $\R^3$ in the
following sense : $D_H(K_{n_k}\cap A,K\cap A)$ tends to $0$ for
every compact set $A$ in $\R^3$. Moreover for all connected
component  $\Omega$ of $\R^3 \backslash K$ and for all compact set
$A$ of $\Omega$, there is a sequence $a_k$ such that
$\{u_{n_k}-a_k\}_{k \in \N}$ converges to $u$ in $L^1(A)$.
 Then, using
\eqref{lim} and Proposition 37.18 of \cite{d}, we know that for
every ball $B \subset \R^3$,
$$\int_{B\backslash K}|\nabla u|^2\leq \liminf_{k \to +\infty}
\int_{B\backslash K_n}|\nabla u_n|^2 \leq\lim_{k\to +\infty}
r\frac{1}{n_k}= 0.$$ Thus $\nabla u=0$ and $u$ is locally
constant. Finally, Theorem  38.3 of \cite{d} says that the limit
$(u,K)$ is a Mumford-Shah minimizer with gauge function $h_l(4r)$.
Since it is true for all  $l$, and that $\sup_lh_l=0$, we can
suppose that $(u,K)$ is a Mumford-Shah minimizer with gauge
function equal to zero, and $u$ is locally constant. But in this
case we know by \cite{d5} that $K$ is a minimal cone of type
$\Pp$, $\Y$ or $\T$, and since for all $n$, $K_n$ is containing
$x$, it is still true for the limit $K$. In addition, there is a
rank $L$ such that for all $k \geq L$ we have $D_H(K\cap
B(0,1),K_{n_k}\cap B(0,1))\leq \frac{\eta_4}{2}$ which is in
contradiction with \eqref{derv} and achieve the proof. \qed

Lemma \ref{compa3} implies the following Theorem.

\begin{theo} \label{fin5} There is some positive constants
$\varepsilon$ and $c<1$ such that the following is true. Let
$(u,K)$ be a Mumford-Shah minimizer in  $\Omega \subset \R^3$ with
gauge function $h$, let $x \in K$ and $r$ be such that
$B(x,r)\subset{\Omega}$ and
$$\omega_2(x,r)+h(r)\leq \varepsilon.$$ Then there is a diffeomorphism
 $\phi$ of class $C^{1,\alpha}$ from  $B(x, cr)$
to its image, and there is a minimal cone $Z$ such that $K \cap
B(x, Cr)=\phi(Z)\cap B(x,cr)$.
\end{theo}

{\bf Proof :} Denote by $\bar \varepsilon$ the constant of Theorem
\ref{fin3}. We apply Lemma \ref{compa3} to $(u,K)$ with
$\eta_4=\bar \varepsilon$. We know that there is a constant $c<1$
and there is a cone $Z$ that contains $x$ such that
$$D_{x,cr}(Z,K)\leq \bar \varepsilon.$$
Dividing if necessary $c$ by $16$ we may assume that the center of
the cone lies in $\frac{1}{4}B(x,cr)$. Thus there is an $y\in
B(x,c\frac{r}{2})$ such that, possibly taking  a smaller
$\varepsilon$,
$$ D_{y,c\frac{r}{2}}(Z,K)+\omega_{2}(y,c\frac{r}{2})+h(r)\leq \bar \varepsilon $$
and then we can apply Theorem \ref{fin3} in $B(y,c\frac{r}{2})$,
and the conclusion follows.  \qed

By the same way of Corollary  \ref{fin4}, in terms of functional
$J$ we have the following statement.

\begin{cor} \label{fin6} There exist some positive constants
$\varepsilon$ and $c$ such that the following is true. Let $g \in
L^\infty$ and $\Omega \subset \R^3$. There is a $\tilde r$
depending only on $\|g\|_{\infty}$, such that for all pair
 $(u,K)\in \mathcal{A}$ that minimizes
$$J(u,K):= \int_{\Omega \backslash K}|\nabla u|^2dx+\int_{\Omega\backslash K}(u-g)^2dx+H^{2}(K),$$
for all $x \in K$ and $r<\tilde r$ such that
$$\omega_2(x,r)\leq \varepsilon$$
 there is a diffeomorphism   $\phi$ of class
 $C^{1,\alpha}$ from $B(x, cr)$ to its image such that $K \cap B(x,
cr)=\phi(Z)\cap B(x,cr)$.
\end{cor}

\bibliographystyle{alpha}
\bibliography{biblio}

ADDRESS :

Antoine LEMENANT \\
e-mail : antoine.lemenant@math.u-psud.fr

Université Paris XI\\
Bureau 15 Bâtiment 430 \\
ORSAY 91400 FRANCE

Tél: 00 33 169157951

\end{document}